\documentclass{aptpub}

\authornames{Fran\c{c}ois Baccelli, Eliza O'Reilly} 
\shorttitle{} 

\usepackage{ dsfont }

\usepackage{xcolor}
\usepackage{changebar}
\newcommand{\R}{\mathds{R}}
\newcommand{\Var}{\mathrm{Var}}
\newif\ifdiff

\ifdiff
\newcommand{\removed}[1]{\cbstart\removedfragile{#1}\cbend{}}
\newcommand{\removedfragile}[1]{{\color{red}{#1}}{}}
\newcommand{\added}[1]{\cbstart\addedfragile{#1}\cbend{}}
\else
\newcommand{\removed}[1]{} 
\newcommand{\added}[1]{#1}

\fi

\allowdisplaybreaks

\begin{document}

\title{Reach of Repulsion for Determinantal Point \\ Processes in High Dimensions} 

\authorone[University of Texas at Austin]{Fran\c{c}ois Baccelli}
\authorone[University of Texas at Austin]{Eliza O'Reilly} 

\addressone{University of Texas at Austin, Department of Mathematics, RLM 8.100,
2515 Speedway Stop C1200
Austin, Texas 78712-1202} 

\begin{abstract}
Goldman [7] proved that the distribution of a stationary determinantal point process (DPP) $\Phi$ can be coupled with its reduced Palm version $\Phi^{0,!}$ such that there exists a point process $\eta$ where $\Phi =  \Phi^{0,!}  \cup \eta$ in distribution and $\Phi^{0,!} \cap \eta = \emptyset$. The points of $\eta$ characterize the repulsive nature of a typical point of $\Phi$. In this paper, the first moment measure of $\eta$ is used to study the repulsive behavior of DPPs in high dimensions. 
It is shown that many families of DPPs have the property that the total number of points in $\eta$ converges in probability to zero as the space dimension $n$ goes to infinity. 
It is also proved that for some DPPs there exists an $R^*$ such that the decay of the first moment measure of $\eta$ is slowest in a small annulus around the sphere of radius $\sqrt{n}R^*$. This $R^*$ can be interpreted as the asymptotic reach of repulsion of the DPP. Examples of classes of DPP models exhibiting this behavior are presented and an application to high dimensional Boolean models is given.
\end{abstract}



\section{Introduction} 
Determinantal point processes (DPPs) are useful models for point patterns where the points exhibit some repulsion from each other, resulting in a more regularly spaced pattern than a Poisson point process. These models originally appeared in random matrix theory and the formalism was introduced by O. Macchi \cite{Macchi} who was motivated by modeling Fermionic particles in quantum mechanics. They have since been used in many applications, such as telecommunication networks, machine learning, ecology, etc. See \cite{Peres}, \cite{Taskar}, \cite{Moller},   \cite{Baccelli_wireless}, and the references therein. This paper describes the repulsive behavior of stationary and isotropic DPPs as the space dimension goes to infinity. 

In the following, a ball with center at the origin and radius $r$ in $\R^n$ is denoted $B_n(r)$. The $\ell^2$ vector norm will be denoted by $|\cdot|$ and the $L^2$-norm on the space $L^2(\R^n)$ by $||\cdot||_2$. Now, consider a sequence of point processes $\Phi_n$ indexed by dimension, each with constant intensity $\rho_n$. If $\rho_n = e^{n \rho}$ and $R_n = \sqrt{n}R$, with $\rho \in \R$ and $R > 0$, then Stirling's formula gives
\begin{equation*}
\text{Vol}(B_n(R_n)) \sim \frac{1}{\sqrt{n\pi}}\left(\frac{2\pi e}{n}\right)^{\frac{n}{2}}R_n^n, \text{ as } n \to \infty.
\end{equation*}
This implies there exists a threshold \added{$R^* = \frac{1}{\sqrt{2\pi e}e^{\rho}}$} such that as $n \rightarrow \infty$,
\begin{equation}\label{Shannon}
\mathds{E}[\Phi_n(B_n(R_n))] \sim e^{n(\rho + \frac{1}{2} \log 2 \pi e + \log R) + o(n)} \rightarrow \begin{cases} 0, & R <  R^* \\ \infty, & R > R^*. \end{cases}
\end{equation}
This justifies the interest in considering this regime where the intensities grow exponentially with dimension and distances grow with the square root of the dimension. This regime also naturally arises in information theory, and following \cite{Venkat}, it will be called the Shannon regime. In this paper, the effect of repulsion in this regime is studied and the range and strength at which DPPs asymptotically exhibit repulsion between points is quantified.

Mention of these issues appear in \cite{Torquato}, where the authors characterize a certain class of DPPs by an effective ``hard-core" diameter $D$ that grows like $\sqrt{n}$, aligning with our observations. They observe that \added{for $r < D$, the number of points in a ball of radius $r$ around a typical point will be zero with probability approaching one, and }for $r > D$, the number of points in a ball of radius $r$ around a typical point is zero with probability approaching zero as dimension $n$ goes to infinity. \added{The behavior for $r < D$ is a result of the natural separation due to dimensionality as exhibited in \eqref{Shannon}.} However, the observation that $D$ is the maximal such separation is due to the \added{$\nu$-weakly sub-Poisson property of DPPs as defined in \cite{Blasz}, and} is a feature of all DPPs, not just those studied in \cite{Torquato}. \added{This behavior is the same as a sequence of Poisson point processes in the same regime, and thus this separation of points in high dimensions is due to dimensionality and not the repulsion of the DPP model.} In this paper, a more precise description of the repulsive behavior in high dimensions is given that is specific to the associated kernel of the DPP.

The measure of repulsiveness used in this paper is a refinement of the global measure of repulsiveness for stationary DPPs described in the on-line supplementary material to \cite{Moller} (see \cite{Moller_Arxiv}). In that work, the authors consider the measure
\begin{equation}\label{gamma}
\gamma := \rho \int \left( 1 - g(x) \right) {\mathrm d}x,
\end{equation}
where $\rho$ is the intensity, and $(x,y) \mapsto g(x-y)$ is the pair correlation function of the point process.
A point process is considered more repulsive the farther $g$ is away from 1; $g \equiv 1$ corresponds to a Poisson point process. \added{As observed in \cite{Kuna}, this measure has the upper bound $\gamma \leq 1$ for all stationary point processes.}

This measure can be refined in order to examine the repulsive effect of a point of the point process across some finite distance. Goldman \cite{Goldman} proved that for a stationary DPP $\Phi$ satisfying certain conditions, there exists a point process $\eta$ such that 
\begin{equation*}
\Phi = \Phi^{0,!} \cup \eta \text{ in distribution, and } \Phi^{0,!} \cap \eta = \emptyset, 
\end{equation*}
where $\Phi^{0,!}$ denotes a point process with the reduced Palm distribution of $\Phi$. Thus, $\eta$ is the set of points that have to be removed from $\Phi$ due to repulsion when a point is ``placed at" the origin.
In the following, the first moment measure of $\eta$ will be used as a measure of the repulsiveness of a DPP $\Phi$, and the repulsive effect of a typical point over a finite distance $R$ is quantified by $\mathds{E}[\eta(B_n(R))]$. Note also that
\begin{equation*}
\mathds{E}[\eta(B_n(R))] = \rho\mathrm{Vol}(B_n(R)) - \mathds{E}[\Phi^{0,!}(B_n(R))] =\rho \left[K_{Poi}(R) - K_{DPP}(R)\right],
\end{equation*}
where $K_{Poi}$ and $K_{DPP}$ are Ripley's K-functions \cite{Ripley} for a Poisson point process and $\Phi$, respectively. Finally, note that the measure of global repulsiveness \eqref{gamma} corresponds to $\eta$ in the sense that $\gamma = \mathds{E}[\eta(\R^n)]$. \added{In recent work \cite{JMEO}, couplings of DPPs and their reduced Palm distributions used to quantify repulsiveness of DPPs are studied further.}

Our main results describe the behavior of the first moment measure of $\eta$ in the Shannon regime. Consider a sequence of stationary DPPs $\{\Phi_n\}$, such that $\Phi_n$ lies in $\R^n$. For each $n$, let $\eta_n$ be the point process such that $\Phi_n = \Phi^{0,!}_n \cup \eta_n$ in distribution and $\Phi_n^{0,!} \cap \eta_n = \emptyset$. One can consider the quantity $\mathds{E}[\eta_n(\R^n)]$ and the probability measure $\frac{\mathds{E}[\eta_n(\cdot)]}{\mathds{E}[\eta_n(\R^n)]}$ on $\R^n$ that is defined to estimate the strength and reach of the repulsiveness of a DPP in any dimension. 

It is often the case that $\mathds{E}[\eta_n(\R^n)] \rightarrow 0$ as $n \rightarrow \infty$. In this case, Markov's inequality and the coupling inequality imply that, in high dimensions, the total variation distance is small between $\Phi_n$ and $\Phi_n^{0,!}$. Indeed, 
\begin{equation}\label{e:coupling_ineq}
|| \Phi_n - \Phi^{0,!}_n ||_{TV} \leq \mathds{P}(\eta_n(\R^n) > 0) \leq \mathds{E}[\eta_n(\R^n)].
\end{equation}
Since $\Phi_n$ and $\Phi_n^{0,!}$ have the same distribution if and only if $\Phi_n$ is Poisson by Slivnyak's theorem \cite{Stoyan}, this says that such DPPs look increasingly like Poisson point processes as the space dimension increases. 

However, the effect of the repulsion can still be observed by examining the probability measure $\frac{\mathds{E}[\eta_n(\cdot)]}{\mathds{E}[\eta_n(\R^n)]}$ on $\R^n$ as seen in Propositions \ref{thresh}, \ref{logconcave}, and \ref{LDP}. Letting $X_n$ be a random vector in $\R^n$ with this probability distribution, it is shown that if $\frac{|X_n|}{\sqrt{n}} \rightarrow R^* \in (0, \infty)$ in probability, then 
\begin{equation*}
\lim_{n \to \infty} \frac{\mathds{E}[\eta_n(B_n(R\sqrt{n}))]}{\mathds{E}[\eta_n(\R^n)]} =\begin{cases} 0, & R < R^* \\ 1, & R > R^*. \end{cases}
\end{equation*}
Here, $R^*$ is interpreted as the {\em asymptotic reach of repulsion} in the Shannon regime for these DPPs. This result implies that in high dimensions, a typical point has its strongest repulsive effect on points that are at a distance of $\sqrt{n}R^*$ away. 

\added{The parametric families of DPP kernels presented in \cite{Lav} and \cite{Moller} provide examples of DPPs exhibiting a reach of repulsion $R^*$ and counterexamples where no finite $R^*$ exists,} as well as computational results on the rates of convergence when a threshold does occur. Four classes of DPPs are studied in Section \ref{examples}: Laguerre-Gaussian DPPs, power exponential DPPs, Bessel-type DPPs, and normal-variance mixture DPPs. For Laguerre-Gaussian DPPs, the sequence ${|X_n|}/{\sqrt{n}}$ satisfies a large deviations principle (established later in Lemma \ref{Laguerre}). As a consequence, the reach of repulsion $R^*$ becomes a phase transition for the exponential rate at which $\mathds{E}[\eta_n\left(B_n(R\sqrt{n})\right)] \rightarrow 0$ as $n \rightarrow \infty$ (established later in Proposition \ref{Laguerre_rate}). Power exponential DPPs are shown to have a finite reach of repulsion in the Shannon regime for certain parameters (established later in Proposition \ref{powerexp_reach}). Bessel-type DPPs are a more repulsive family that does not exhibit an $R^*$ (established later in Proposition \ref{bessel}). Finally, normal-variance mixture DPPs provide additional examples of DPPs that exhibit an $R^*$, including the Cauchy and Whittle-Mat\'ern models (established later in Propositions \ref{Cauchy} and \ref{Whittle}). 

An application of these results is presented in Section \ref{applications}. It can be shown that some threshold results in \cite{Venkat} for Poisson Boolean models can be extended to generalized Laguerre-Gaussian DPP Boolean models in the Shannon regime using the rates of convergence computed for these DPPs. Finally, concluding remarks and open questions are stated in Section \ref{conclusion}.


\section{Preliminaries}\label{DPPs}


Determinantal point processes are characterized by an integral operator $\mathcal{K}$ with kernel $K$, and can be defined in terms of their joint intensities, also known as correlation functions (\cite{Peres}, \cite{Moller}). 

\begin{defn} A simple, locally finite, spatial point process $\Phi$ on $\mathds{R}^n$ is a determinantal point process with kernel $K: \R^n \times \R^n \rightarrow \mathds{R}$ ($\Phi \sim DPP(K)$) if its joint intensities exist for all order $k$ and satisfy
\begin{align*}
\rho^{(k)}(x_1, \dotsc , x_k) = \det(K(x_i, x_j))_{1 \leq i,j, \leq k}, \qquad k = 1,2, \dotsc .
\end{align*}
\end{defn}

Note that the intensity function of $\Phi$ is given by $\rho(x) = K(x,x)$. The degenerate case where $K(x,y) = \delta_{\{x = y\}}$ coincides with a Poisson point process with unit intensity. 

The following conditions on $K$ are imposed to ensure $\Phi \sim DPP(K)$ is well-defined. 
Let $K: \R^n \times \R^n \rightarrow \mathds{R}$ be a continuous 
kernel and assume $K$ is symmetric, i.e., $K(x,y) = K(y,x)$. The kernel $K$ then defines a self-adjoint integral operator $\mathcal{K}$ on $L^2(\R^n)$ given by $\mathcal{K}f(x) = \int K(x,y )f(y) {\mathrm d}y$. For any compact set $S \subset \R^n$, the restricted operator $\mathcal{K}_S$ given by
\begin{align*}
\mathcal{K}_Sf(x) = \int_S K(x,y)f(y) {\mathrm d}y, \qquad \text{ } x \in S,
\end{align*}
is a compact operator. By the spectral theory for self adjoint compact operators, the spectrum of $\mathcal{K}_S$ consists solely of countably many eigenvalues $\{\lambda_k^S\}_{k \in \mathds{N}}$ with an accumulation point only possible at zero. See \cite{Rudin} for more on compact operators.
These conditions imply that for any compact $S \subset \R^n$, the kernel $K$ restricted to $S \times S$ has a spectral representation 
\begin{align*}
K(x,y) = \sum_{k=1}^{\infty} \lambda_k^S \phi_k^S(x) \overline{\phi_k^S(y)}, \qquad (x,y) \in S \times S,
\end{align*}
where $\{\phi_k^S\}_{k \in \mathds{N}}$ are the eigenvectors of $\mathcal{K}_S$, and form an orthonormal basis of $L^2(S)$. 

\begin{thm}\label{specbound}(\added{Macchi \cite{Macchi}}) Under the conditions given above, a kernel $K$ defines a determinantal process on $\mathds{R}^n$ if and only if the spectrum of $\mathcal{K}$ is contained in $[0,1]$.
\end{thm} 

If $K(x,y) = K_0(x-y)$, then $\Phi \sim DPP(K)$ is stationary. In this case, the operator $\mathcal{K}$ is the convolution operator $\mathcal{K}(f) = K_0 \star f$ on $L^2(\R^n)$. The intensity function $\rho(x)$ is then constant and satisfies $\rho = K_0(0)$. 
For these stationary DPPs, there is a simple spectral condition for existence.

\begin{thm}\label{exist} (Theorem 2.3 in \cite{Moller}) Assume $K_0$ is a symmetric continuous real-valued function in $L^2(\mathds{R}^n)$. Let $K(x,y) = K_0(x-y)$. Then DPP($K$) exists if and only if $0 \leq \hat{K}_0 \leq 1$, where $\hat{K}_0$ denotes the Fourier transform of $K_0$. 
\end{thm} 

For the rest of this paper, when it is stated that $\Phi \sim DPP(K)$ is stationary, it is assumed that $K(x,y) = K_0(x - y)$ for a real-valued $K_0 \in L^2(\R^n)$, and $K$ will be used to mean $K_0$. There exist stationary DPPs with kernels that are not of this form (see \cite[4.3.7]{Peres}), but they are complex-valued and not considered here. In addition, when it is stated that $\Phi$ is isotropic, it is meant that $K_0(x) = R_0(|x|)$ and the distribution of $\Phi$ is thus invariant under rotations about the origin in $\R^n$.

The reduced Palm distribution of a stationary point process $\Phi$ can be interpreted as the distribution of $\Phi$ conditioned on there being a point at the origin with the point at the origin removed (see \cite[Chapter 4]{Stoyan}) and will be denoted by \added{$\mathds{P}^{0,!}$. A point process with the Palm distribution $\mathds{P}^{0,!}$ of $\Phi$ will be denoted $\Phi^{0,!}$.}  The following theorem is a special case of a useful result about the Palm distribution of DPPs. 

\begin{thm}\label{palm} (Theorem 6.5 in \cite{Shirai}). Let $\Phi \sim DPP(K)$ in $\R^n$ be stationary with intensity $\rho = K(0) > 0$. Then $\Phi^{0,!}$ is a DPP with associated kernel
\begin{center}
$K^!_{0} (x,y) = \frac{1}{K(0)} \det\left(\begin{array}{cc} K(x - y) & K(x ) \\ K(y) & K(0) \end{array}\right) = K(x-y) - \frac{1}{\rho}K(x)K(y)$.
\end{center} 
\end{thm}

The nearest neighbor function of a stationary point process $\Phi$ in $\R^n$ is defined as
\begin{equation}\label{neighbor}
D(r) := \mathds{P}^{0,!}(\Phi(B_n(r)) > 0).
\end{equation}
If $\Phi$ is Poisson, Slivnyak's theorem gives that $D(r) = 1 - e^{-\mathds{E}\Phi(B_n(r))}$. For $\Phi \sim DPP(K)$, Theorem \ref{palm} implies that $D(r) = \mathds{P}(\Phi^{0,!}(B_n(r)) > 0)$, with $\Phi^{0,!} \sim DPP(K^!_0)$.

As mentioned in the introduction, Goldman \cite{Goldman} proved the following result. 
\begin{thm}\label{eta} (Theorem 7 in \cite{Goldman}) 
Let $\Phi \sim DPP(K)$, where $K$ is continuous, and the spectrum of the integral operator $\mathcal{K}$ with kernel $K$ is contained in $[0,1)$. Then, there exists a point process $\eta$ such that
\begin{align*}
\Phi = \Phi^{0,!} \cup \eta \text{ in distribution, and } \Phi^{0,!} \cap \eta = \emptyset.
\end{align*}
\end{thm}

This theorem says that a point process with the distribution of $\Phi^{0,!}$ can be obtained from $\Phi$ by removing a subset of points $\eta$. This is a striking result, since the procedure does not include shifting any of the remaining points. The points in $\eta$ characterize the repulsive nature of the DPP $\Phi$, since these are the points that are ``pushed out" by the point at zero under the reduced Palm distribution. It also makes sense to compare the repulsiveness of DPPs using $\eta$. For two stationary DPPs $\Phi_1$ and $\Phi_2$ with the same intensity, $\Phi_1$ is defined to be more repulsive than $\Phi_2$ if $\mathds{E}[\eta_1(\R^n)] > \mathds{E}[\eta_2(\R^n)]$. This corresponds to the definition in \cite{Moller} using the measure $\gamma$ defined in (\ref{gamma}).
Note that the assumptions for Theorem \ref{eta} excludes the interesting case where $\mathcal{K}$ has an eigenvalue of 1, which corresponds to when $\hat{K}(x)$ attains a value of one for some $x$.

\section{Main Results}\label{main}

When considering the reach of repulsion of a DPP, it is natural to first consider the nearest neighbor function (\ref{neighbor}). The following threshold behavior was observed for stationary DPPs in \cite{Torquato}. It is stated here for a sequence of DPPs in the Shannon regime. For each $n$, let $\Phi_n \sim DPP(K_n)$ in $\R^n$ be stationary with intensity $K_n(0) = e^{n\rho}$ for some $\rho \in \R$. Then, for $\tilde{R} :=  (2\pi e)^{-\frac{1}{2}} e^{-\rho}$,
\begin{align}\label{near_thresh}
\lim_{n \rightarrow \infty} \mathds{P}(\Phi^{0,!}_n(B_n(\sqrt{n}R)) > 0) = \begin{cases} 0, & R < \tilde{R} \\ 1, & R > \tilde{R} .\end{cases}
\end{align}
A proof of this fact is given in Appendix \ref{near_thresh_app}.

This shows there is a separation of points as dimension tends to infinity for any stationary DPP. However, the same threshold behavior occurs if the elements of the sequence $\{\Phi_n\}$ are stationary Poisson point processes, as a consequence of \added{\eqref{Shannon}}. This observation shows that \added{this separation is due purely to dimensionality and is not a result of the repulsiveness of DPPs.}

\added{The point process $\eta_n$ as defined in Theorem \ref{eta} gives an alternative characterization of the repulsiveness of a DPP and can measure some consequence of repulsiveness in high dimensions that depends on the determinantal structure.}

\begin{lem}\label{firstmoment} \added{Let $\Phi_n \sim DPP(K_n)$ in $\R^n$ be stationary and assume $0 \leq \hat{K}_n < 1$. Let $\eta_n$ be the point process given in Theorem \ref{eta} and define the random vector $X_n$ in $\R^n$ with probability density $\frac{K_n(x)^2}{||K_n||_2^2}$. Then,
\[ \mathds{P}(X_n \in B) = \frac{\mathds{E}[\eta_n(B)] }{\mathds{E}[\eta_n(\R^n)] }, \qquad B \in \mathcal{B}(\R^n).\]}
\end{lem}

The following result shows that under certain limit conditions on the kernels of a sequence of DPPs, the repulsiveness measured by the first moment measure of $\eta_n$ is concentrated at a distance of $\sqrt{n}R^*$ for some $R^* \in (0, \infty)$ as $n$ goes to infinity. 

\begin{prop}\label{thresh}
For each $n$, let $\Phi_n \sim DPP(K_n)$ be a stationary and isotropic DPP in $\R^n$, and assume $0 \leq \hat{K}_n < 1$. Let $X_n$ be a random vector in $\R^n$ with probability density $\frac{K_n(x)^2}{||K_n||_2^2}$. Assume that \added{as $n \rightarrow \infty$,
\begin{align}\label{conv}
\frac{|X_n|}{\sqrt{n}} \rightarrow R^* \text{ in probability.}
\end{align}}
Then, 
\begin{align}\label{threshold}
\lim_{n \rightarrow \infty} \frac{\mathds{E}[\eta_n(B(\sqrt{n}R))]}{\mathds{E}[\eta_n(\R^n)]} = \begin{cases} 0, & R < R^* \\ 1, & R > R^*. \end{cases} 
\end{align}
\end{prop}

\added{\begin{rem}\label{chebychev}
One way to show \eqref{conv} is to show that \[\lim_{n \rightarrow \infty} \frac{\Var(|X_n|^2)}{n^2} = 0 \text{ and } \lim_{n \rightarrow \infty} \left( \frac{\mathds{E}[|X_n|^2]}{n} \right)^{1/2} = R^* \in (0,\infty),\] and then apply Chebychev's inequality. 
\end{rem}}

\begin{rem}\label{thinshell}
For general vectors $X_n$ in $\R^n$, the concentration of \added{$|X_n|$ for large $n$ has been well-studied (see \cite{Guedon}, \cite{Milman}, \cite{Klartag}). Indeed, in \cite[Proposition 3]{Guedon}, it is proved that $X_n$ is concentrated in a ``thin shell", i.e., there exists a sequence $\{\varepsilon_n\}$ such that $\varepsilon_n \rightarrow 0$ as $n \rightarrow \infty$ and for each $n$,
\begin{align}\label{thinshellinequality}
\mathds{P}\left(\bigg|\frac{|X_n|}{\mathds{E}[|X_n|^2]^{\frac{1}{2}}} - 1 \bigg| \geq \varepsilon_n \right) \leq \varepsilon_n,
\end{align}}
if and only if $|X_n|$ has a finite $r$th moment for $r > 2$, and for some $2 < p < r$,
\begin{align*}
\bigg|\frac{\mathds{E}[|X_n|^p]^{1/p}}{\mathds{E}[|X_n|^2]^{1/2}} - 1 \bigg| \rightarrow 0 \text{ as } n \rightarrow \infty.
\end{align*}
\end{rem} 

For random vectors with log-concave distributions, the deviation estimate can be improved from the estimate obtained through Chebychev's inequality (see \added{Remark \ref{chebychev}}). The best known estimate is given by the following theorem in \cite{Milman}. 

\begin{thm}\label{thinshellthm} (Gu\'{e}don and Milman \cite{Milman}) Let $X$ denote \added{a random vector in $\R^n$ such that $\mathds{E}X = 0$ and $\mathds{E} (X \otimes X) = I_n$.} Assume $X$ has a log-concave density. Then, for some $C > 0$ and $c > 0$,
\begin{align*}
\mathds{P}\left( \bigg| \frac{ |X|}{\sqrt{n}} - 1  \bigg| \geq t \right) \leq C \exp\left(-c \added{\sqrt{n}} \min(t^3, t) \right). 
\end{align*}
\end{thm}

This gives the following result. 

\begin{prop}\label{logconcave}
For each $n$, let $\Phi_n \sim DPP(K_n)$ be a stationary and isotropic DPP in $\R^n$, and assume $0 \leq \hat{K}_n < 1$. Let $X_n$ be a random vector with density $\frac{K_n(x)^2}{||K_n||_2^2}$\added{ and let $\sigma_n^2 = \mathds{E}|X_n|^2$.} If $K_n^2$ is log-concave for all $n$, 
then there exist positive constants $C, c$ such that for all $\delta \in (0, 1)$,
\begin{equation*}
\frac{\mathds{E}[\eta_n(B_n(\sigma_n(1-\delta)))]}{\mathds{E}[\eta_n(\R^n)]} \leq C e^{-c \added{\sqrt{n}}\delta^3 },
\end{equation*}
and for all $\delta > 0$,
\begin{equation*}
 \frac{\mathds{E}[\eta_n(\R^n\backslash B_n(\sigma_n(1 + \delta)))]}{\mathds{E}[\eta_n(\R^n)]} \leq C e^{-c \added{\sqrt{n}} \min(\delta^3, \delta)}.
\end{equation*} 
If, in addition, 
\begin{align}\label{e:meanconv}
\lim_{n \rightarrow \infty} \frac{\sigma_n}{\sqrt{n}} = R^* \in (0, \infty),
\end{align}
then for this $R^*$, the threshold \eqref{threshold} occurs, 
and for all $R < R^*$, there exists a constant $C(R) > 0$ such that
\begin{align*}
\liminf_{n \rightarrow \infty} - \frac{1}{\added{\sqrt{n}}} \ln \frac{\mathds{E}[\eta_n(B_n(\sqrt{n}R))]}{\mathds{E}[\eta_n(\R^n)]} \geq C(R).
\end{align*}
\end{prop}

\begin{rem}
The last conclusion of Proposition \ref{logconcave} about the rate also holds for $R > R^*$ if $B_n(\sqrt{n}R)$ is replaced by $\R^n \backslash B_n(\sqrt{n}R)$.
\end{rem}

The assumption of large deviation principle (LDP) concentration leads to an estimate of the exponential rate of convergence with speed $n$ and an exact computation of the reach of repulsion $R^*$.
\begin{prop}\label{LDP}
For each $n$, let $\Phi_n \sim DPP(K_n)$ be a stationary and isotropic DPP in $\R^n$, and assume $0 \leq \hat{K}_n < 1$. Let $X_n$ be a random vector with density $\frac{K_n(x)^2}{||K_n||_2^2}$ and suppose $\frac{|X_n|}{\sqrt{n}}$ satisfies a LDP with strictly convex rate function $I$. Then, for $R^*$ such that $I(R^*) = 0$, the threshold \eqref{threshold} occurs.
Also, for $R < R^*$,
\begin{equation*}
- \inf_{r < R} I(r) \leq \liminf_{n \rightarrow \infty}  \frac{1}{n} \ln \frac{\mathds{E}[\eta_n(B_n(\sqrt{n}R))]}{\mathds{E}[\eta_n(\R^n)]} \leq \limsup_{n \rightarrow \infty} \frac{1}{n} \ln \frac{\mathds{E}[\eta_n(B_n(\sqrt{n}R))]}{\mathds{E}[\eta_n(\R^n)]} \leq  - \inf_{r \leq R} I(r), 
\end{equation*}
and if the rate function $I$ is continuous at $R$,
\begin{equation*}
\lim_{n \rightarrow \infty} - \frac{1}{n} \ln \frac{\mathds{E}[\eta_n(B_n(\sqrt{n}R))]}{\mathds{E}[\eta_n(\R^n)]} = I(R).
\end{equation*}
\end{prop}

\begin{rem}
The second conclusion of Proposition \ref{LDP} about the rate also holds for $R > R^*$ if $B_n(\sqrt{n}R)$ is replaced by $\R^n \backslash B_n(\sqrt{n}R)$.
\end{rem}

If a sequence of DPPs in increasing dimensions exhibits a reach of repulsion $R^*$, this says that the points of $\eta_n$ are most likely to be near distance $\sqrt{n}R^*$ away from the origin in high dimensions. If $R^*$ is less than $\tilde{R}$ from \eqref{near_thresh}, points are most likely to be removed at a distance where points of $\Phi_n$ appear with probability decreasing to zero as $n$ increases due to dimensionality. If $R^*$ can reach past $\tilde{R}$, the points ``pushed out" by repulsion are most likely to lie at a distance where points of $\Phi_n$ appear with high probability. Thus it is of interest to check whether there exist DPP models such that $R^*$ is greater than or equal to $\tilde{R}$, i.e., if $\mathds{P}(\Phi^{0,!}_n(B_n(\sqrt{n}R^*)) = 0) \rightarrow 0$ as $n \rightarrow \infty$. In Sections \ref{laguerre_section} and \ref{powerexp} examples of DPP models with this reach are provided. 

The above results have strong assumptions, and open up additional questions. The first question is whether the points of $\eta_n$ tend to lie at distances scaling with $\sqrt{n}$, i.e., is the Shannon regime the right one to examine the repulsiveness between points of a family of DPPs in high dimensions? By the radial symmetry of the density of each $X_n$, the coordinates $\{X_{n,k}\}_{k=1}^n$ are identically distributed, and the sequence $|X_n|^2$ is the sequence of row sums of a triangular array of random variables with identically distributed rows. If the coordinate distributions depend on dimension in such a way that $\mathds{E}\left(|X_n|^2\right) \neq O(n)$, then a different scaling is needed. 

\section{Examples}\label{examples}

\added{In the following, specific families that were presented in \cite{Lav} and \cite{Moller} are examined that illustrate both examples of DPP models satisfying the above results, as well as examples that do not. These examples provide a window into the wide scope of repulsive behavior that can be described using this framework. }

The first task will be to determine the behavior of $\mathds{E}[\eta_n(\R^n)]$ as $n$ increases. For each of the examples provided in this section, $\lim_{n \rightarrow \infty} \mathds{E}[\eta_n(\R^n)] = 0$, but each class exhibits this convergence at different speeds. Then the goal is to determine if the DPP models satisfy the conditions of Propositions \ref{thresh}, \ref{logconcave}, or \ref{LDP}. 

\subsection{Laguerre-Gaussian Models}\label{laguerre_section}

\added{For each $n$, let $\Phi_n \sim DPP(K_n)$ in $\R^n$ be a Laguerre-Gaussian DPP as described in \cite{Lav} with intensity $K_n(0) = e^{n\rho}$, i.e., for some $m \in \mathds{N}$, $\alpha \in \R^+$, let}
\begin{align}\label{e:Lag_K}
K\added{_n}(x) = \frac{e^{n\rho}}{\binom{ m - 1 + \frac{n}{2}}{m - 1}} L_{m-1}^{n/2}\left(\frac{1}{m} \left|\frac{x}{\alpha}\right|^2\right) e^{-\frac{|x/\alpha|^2}{m}},\qquad x \in \R^n,
\end{align}
\added{where $L_m^{\beta}(r) = \sum _{k=0}^{m} \binom{m + \beta }{ m - k} \frac{(-r)^k}{k!}$, for all $r \in \R$, denote the Laguerre polynomials.}
\added{From \cite{Lav}, the} condition $0 \leq \hat{K}\added{_n} < 1$ translates to a bound on $\alpha\added{_n}$,
\added{\begin{align}\label{e:alpha_bnd} \alpha <  \frac{1}{e^\rho (m\pi)^{1/2}} \binom{m - 1 + n/2 }{m - 1}^{\frac{1}{n}}.
\end{align}}
Direct calculations give that the global measure of repulsiveness is
\begin{align}\label{Lag_norm}
&\mathds{E}[\eta_n(\R^n)] = \nonumber\\
&\frac{e^{n\rho}\alpha_n^n}{\binom{m - 1 + \frac{n}{2}}{m - 1 }^2} \left(\frac{m\pi}{2}\right)^{\frac{n}{2}} \sum_{k, j=0}^{m-1}\binom{m - 1 + \frac{n}{2}}{m -1 - k }\binom{m - 1 + \frac{n}{2} }{m - 1 - j } \frac{(-1)^{k+j}}{k!j!}   \frac{\Gamma\left(\frac{n}{2} + k + j \right)}{2^{k + j} \Gamma\left(\frac{n}{2}\right)}.
\end{align}
By \added{\eqref{e:alpha_bnd},} $\mathds{E}[\eta(\R^n)]  <  2^{-\frac{n}{2}} f(n,m)$,
where 
\[f(n,m) =  \sum_{k, j=0}^{m-1} \frac{  \binom{ m - 1 + n/2}{ m -1- k}\binom{m - 1 + n/2 }{m -1- j }}{\binom{ m - 1 + n/2 }{m - 1 }}  \frac{(-1)^{k+j}}{k!j!}   \frac{\Gamma\left(\frac{n}{2} + k + j \right)}{2^{k + j} \Gamma\left(\frac{n}{2}\right)} = \added{O(n^{m-1})}.
\]
\added{It follows from \cite[(5.7)]{Lav} that for fixed $n$, $\lim_{m \rightarrow \infty} 2^{-\frac{n}{2}}f(n,m) = 1$, and as $\alpha \rightarrow 0$, $K_n$ approaches the Poisson kernel.} Thus, this class of DPPs covers a wide range of repulsiveness for fixed dimension $n$. However, for any fixed $m$, the dominant behavior as $n$ increases is $2^{-\frac{n}{2}}$.

Since $\binom{m - 1 + n/2 }{ m - 1 }^{\frac{1}{n}}$ decreases to one as $n$ goes to infinity, \added{a sufficient condition for \eqref{e:alpha_bnd} to hold for all $n$ is }$0 < \alpha <  e^{-\rho} (m\pi)^{-\frac{1}{2}}$. Note that this scaling for the intensity is the right one for observing interactions between the parameters of the model because it provides a trade-off between how large the parameter $\alpha$ can be and the magnitude of $\rho$. If the intensity did not grow as quickly with dimension, the upper bound on $\alpha$ would depend less and less on changes in $\rho$ as dimension increased, and if the intensity grew more quickly, the upper bound for $\alpha$ would tend to zero as $n$ goes to infinity.

Proposition \ref{LDP} holds for this sequence of DPPs. Indeed, the next lemma shows that the sequence of $\R^+$-valued random variables $\frac{|X_n|}{\sqrt{n}}$ satisfies a LDP.

\begin{lem}\label{Laguerre} Fix $m \in \mathds{N}$, $\rho \in \R$, and let $\alpha \in (0, e^{-\rho}(m\pi)^{-1/2})$. For each $n$, let $X_n$ be a random vector in $\R^n$ with probability density $\frac{K_n(x)^2}{||K_n||^2_{2}}$, where \added{$K_n$ is given by \eqref{e:Lag_K}.}
Then, the sequence $\{\frac{|X_n|}{\sqrt{n}}\}_n$ satisfies an LDP with rate function 
\begin{align*}
\Lambda^*(x) = \frac{2x^2}{\alpha^2 m} - \frac{1}{2} +  \frac{1}{2} \log \left(\frac{\alpha^2 m}{4x^2}\right).
\end{align*}
\end{lem}

Using this lemma, Proposition \ref{LDP} implies that an $R^*$ exists, and the exponential rates can be determined.
In addition, using (\ref{Lag_norm}), the exponential rate of decay of $\mathds{E}[\eta_n(B_n(\sqrt{n}R))]$ can be computed.

\begin{prop}\label{Laguerre_rate} Fix $m \in \mathds{N}$, $\rho \in \R$, and let $\alpha \in (0, e^{-\rho}(m\pi)^{-1/2})$. For each $n$, let $\Phi_n \sim DPP(K_n)$ where $K_n$ is given by (\ref{e:Lag_K}). 
Then, for $R^* := \sqrt{m}\frac{\alpha}{2}$,
\begin{align*}
\lim_{n \rightarrow \infty} - \frac{1}{n} \log \mathds{E}[\eta_n(B_n(\sqrt{n}R))] = \begin{cases} - \rho - \frac{1}{2} \log 2\pi e  + \frac{2R^2}{\alpha^2m} - \log R, & 0 < R < R^* \\ - \rho - \log \alpha - \frac{1}{2} \log \frac{m\pi}{2} , & R > R^*.
\end{cases}
\end{align*}
\end{prop}

The rate decays as $R$ increases to $R^* := \sqrt{m}\frac{\alpha}{2}$ and then for $R > R^*$, the rate no longer depends on $R$. This coincides with our interpretation of $R^*$ as the asymptotic reach of repulsion of the sequence of DPPs. 

For a fixed $\alpha$, a larger $m$ will give farther reach, and for a fixed $m$, a larger $\alpha$ will provide a farther reach. However, by the bound $\alpha <\frac{1}{e^{\rho} (m\pi)^{1/2}}$, the following upper bound on the reach holds uniformly for all $m$:
\begin{align*}
R^* := \sqrt{m}\frac{\alpha}{2} < \frac{1}{2 e^{\rho} \pi^{1/2}}.
\end{align*}
Note that the larger $\rho$ is, the smaller the upper bound on $R^*$ can be. This follows from the relationship between $\alpha$ and $\rho$: the higher the intensity, the smaller $\alpha$ must be for the DPP to exist. Since a larger $\alpha$ implies a larger values of $\mathds{E}[\eta_n(\R^n)]$, the parameter $\alpha$ is associated with the strength of the repulsiveness. The relationship with $\rho$ showcases the following tradeoff \added{observed in \cite{Moller}}: the higher the intensity of the DPP, the less repulsive it can be.

As mentioned in the previous section, it is of interest to know whether there is a range of parameters such that $R^*$ is greater than $\tilde{R}$, the threshold for the convergence of the nearest-neighbor function of $\Phi$ (\ref{near_thresh}). For Laguerre-Gaussian models, $R^* := \frac{\sqrt{m}\alpha}{2}$ is larger than $\tilde{R}$ and \added{$\alpha$ satisfies the condition of Lemma \ref{Laguerre} if}
\begin{align*}
\left(\frac{2}{e}\right)^{1/2} < e^{\rho}  \sqrt{m\pi} \alpha < 1. 
\end{align*}
Since the lower bound is strictly less than one, there is a non-empty range for $\alpha$ such that the reach of repulsion reaches \added{past $\tilde{R}$.}

\subsection{Power Exponential Spectral Models}\label{powerexp}

The power exponential spectral models, introduced in \cite{Moller}, are defined through the Fourier transform of the kernel. For almost all of these models, there is no closed form for the kernel $K$. Using properties of the Fourier transform, a similar analysis of the repulsive behavior can still be performed.

\added{For each $n$, let $\Phi_n \sim DPP(K_n)$ be a power exponential DPP with intensity $K_n(0) = e^{n\rho}$ and parapmeters $\nu > 0$ and $\alpha_n > 0$, i.e., let}
\begin{align}\label{power_kernel}
\hat{K}_n(x) = e^{n\rho} \frac{\Gamma(\frac{n}{2} + 1) \alpha_n^n}{\pi^{n/2}\Gamma(\frac{n}{\nu} + 1)} e^{-|\alpha_n x|^\nu},\qquad x \in \R^n.
\end{align}
When $\nu = 2$, a closed form expression for $K_n$ exists and is called the Gaussian kernel. 
\added{The condition $0 \leq \hat{K}_n < 1$ implies the following upper bound on $\alpha_n$:
\begin{align}\label{alpha_bnd}
\alpha_n < \frac{\Gamma(\frac{n}{\nu} + 1)^{\frac{1}{n}}\pi^{1/2}}{e^\rho  \Gamma\left(\frac{n}{2} + 1\right)^{\frac{1}{n}}}, 
\end{align}}
and the asymptotic expansion for the upper bound on $\alpha_n$ as $n \rightarrow \infty$ is
\begin{align*}
\left(  \tfrac{\Gamma(\frac{n}{\nu} + 1)\pi^{n/2}}{e^{n\rho} \Gamma(\frac{n}{2} + 1)} \right)^{1/n} 
\sim \left(  \tfrac{\sqrt{\frac{2\pi n}{\nu}} \left(\frac{n}{\nu e}\right)^{n/\nu} \pi^{n/2}}{e^{n\rho} \sqrt{\frac{2\pi n}{2}} \left(\frac{n}{2e}\right)^{n/2} }\right)^{1/n} 
\sim e^{-\rho}n^{\frac{1}{\nu} - \frac{1}{2}}  \frac{(2\pi e)^{1/2}}{(\nu e)^{1/\nu}} = O(n^{\frac{1}{\nu} - \frac{1}{2}}).
\end{align*}
\added{By Parseval's theorem and a change of variables, 
\begin{align}\label{powerexp_global}
\mathds{E}[\eta_n(\R^n)] &= \frac{1}{e^{n\rho}} ||K_n||_2^2 = \frac{1}{e^{n\rho}}||\hat{K}_n||_2^2 = \frac{1}{e^{n\rho}}\left( e^{n\rho} \frac{\Gamma\left(\frac{n}{2} + 1\right) \alpha_n^n}{\pi^{n/2}\Gamma\left(\frac{n}{\nu} + 1\right)} \right)^2 \int_{\R^n} e^{-2|\alpha_n x|^{\nu}}{\mathrm d}x \nonumber \\
&=  e^{n\rho} \left( \frac{\Gamma(\frac{n}{2} + 1) \alpha_n^n}{\pi^{n/2}\Gamma(\frac{n}{\nu} + 1)} \right)^2 \frac{n\pi^{n/2}}{\Gamma(\tfrac{n}{2} + 1)} \int_0^{\infty} r^{n-1}e^{-2(\alpha r)^{\nu}} {\mathrm d}r \nonumber \\
&= e^{n\rho} \frac{\Gamma(\frac{n}{2} + 1) \alpha_n^{2n}}{\pi^{\frac{n}{2}}\Gamma(\frac{n}{\nu} + 1)^2}  \frac{n}{2^{\frac{n}{\nu}} \nu \alpha_n^n} \int_0^{\infty} t^{\frac{n}{\nu} - 1} e^{-t} {\mathrm d}t  
 = 2^{-\frac{n}{\nu}} \alpha_n^{n}\frac{e^{n\rho} \Gamma(\frac{n}{2} + 1) }{\pi^{\frac{n}{2}}\Gamma(\frac{n}{\nu} + 1)}. 
\end{align}
By the bound on $\alpha_n$ \eqref{alpha_bnd}, 
\begin{equation*}
\mathds{E}[\eta_n(\R^n)] < 2^{-\frac{n}{\nu}}.
\end{equation*}
For fixed dimension $n$, the global measure of repulsion approaches its upper bound of one for large $\nu$. Thus, this class covers a wide range of repulsiveness similar to the Laguerre-Gaussian DPPs. However, for fixed $\nu$, the measure decays exponentially as $n$ goes to infinity. Note that for $\nu > 2$, the rate is smaller than for the Laguerre-Gaussian models, i.e., the decay is slower. }

The following results show that if the parameters $\alpha_n$ grow appropriately with $n$, this sequence satisfies the assumptions of Proposition \ref{thresh}.

\begin{lem}\label{scale} For each $n$, let $X_n$ be a vector in $\R^n$  with density $\frac{K_n^2}{||K_n||_2^2}$ such that $\hat{K}_n$ is given by (\ref{power_kernel}).
Assume $\alpha_n \sim \alpha n^{\frac{1}{\nu} - \frac{1}{2}}$ as $n \rightarrow \infty$ for $\alpha \in (0,\infty)$, and $\alpha_n <  \left(  \tfrac{\Gamma(\frac{n}{\nu} + 1)\pi^{n/2}}{e^{n\rho} \Gamma(\frac{n}{2} + 1)} \right)^{1/n} $ for all $n$. Then, as $n \rightarrow \infty$,
\added{\[ \frac{|X_n|}{\sqrt{n}} \rightarrow \alpha \frac{(2\nu)^{1/\nu}}{4\pi} \text{ in probability}.\]}
\end{lem}

Now, applying Proposition \ref{thresh}, the following holds for a sequence of power exponential DPPs in the Shannon regime. 

\begin{prop}\label{powerexp_reach}
For each $n$, let $\Phi_n \sim DPP(K_n)$ where $\hat{K}_n$ satisfies the assumptions in Lemma \ref{scale}. Then, \added{for $R^* := \alpha \frac{(2\nu)^{1/\nu}}{4\pi}$,
\begin{equation*}
\lim_{n \rightarrow \infty} \frac{\mathds{E}[\eta_n(B_n(\sqrt{n}R))]}{\mathds{E}[\eta_n(\R^n)]} = \begin{cases} 0, & R < R^*, \\ 1, & R > R^*. \end{cases}
\end{equation*}}
\end{prop}

For $\nu > 1$, the reach of repulsion $R^*$ for the power exponential models can also reach past the nearest neighbor threshold $\tilde{R}$. Indeed, for $\alpha_n \sim \alpha n^{\frac{1}{\nu} - \frac{1}{2}}$, $R^*:= \alpha \frac{(2\nu)^{1/\nu}}{4\pi}$ satisfies $\mathds{P}[\Phi_n(B_n(0, \sqrt{n}R^*)) = 0] \rightarrow 0$ as $n \rightarrow \infty$ if 
\begin{equation*}
\alpha \frac{(2\nu)^{1/\nu}}{4\pi} >  \frac{1}{\sqrt{2\pi e} e^{\rho}}.
\end{equation*}
By the asymptotic formula (\ref{alpha_bnd}) for the upper bound of $\alpha_n$,
$\alpha < \frac{\sqrt{2\pi e}}{e^{\rho} (\nu e)^{1/\nu}}$.
Thus, $R^*$ reaches \added{past $\tilde{R}$} when $\alpha_n \sim \alpha n^{\frac{1}{\nu} - \frac{1}{2}}$ and
\begin{equation*}
\frac{4\pi}{(2\nu)^{1/\nu}e^{\rho}\sqrt{2\pi e}} < \alpha < \frac{\sqrt{2\pi e}}{e^{\rho} (\nu e)^{1/\nu}}.
\end{equation*}
The interval is non-empty since the upper bound is strictly greater than the lower bound for $\nu > 1$.

\subsection{Bessel-type Models}

Another class of DPP models presented in \cite{Lav} is the Bessel-type. This class is more repulsive than the previous two families of models. It is shown that while the Shannon regime is the right scaling to examine the repulsiveness of this class in high dimensions, a sequence of these DPPs does not satisfy the conditions of Proposition \ref{thresh}.

\added{For each $n$, let $\Phi_n \sim DPP(K_n)$ be a Bessel-type DPP with parameters $\sigma \geq 0$, $\alpha > 0$, and intensity $K_n(0) = e^{n\rho}$, for $\rho \in \R$. That is, let}
\begin{align}\label{e:Bessel_K}
K_n(x) = e^{n\rho} 2^{(\sigma + n)/2} \Gamma\left( \frac{\sigma + n + 2}{2}\right) \frac{J_{(\sigma + n)/2}(2|x/\alpha|\sqrt{(\sigma + n)/2})}{(2|x/\alpha|\sqrt{(\sigma + n)/2})^{(\sigma + n)/2)}}.
\end{align}
\added{From \cite{Lav},} the bound $0 \leq \hat{K}_n < 1$ implies that
\begin{align}\label{a_bnd_bessel}
\alpha_n^n < \frac{(\sigma + n)^{n/2}\Gamma\left(\frac{\sigma}{2} + 1\right)}{e^{n\rho} (2\pi)^{n/2} \Gamma\left(\frac{\sigma + n}{2} + 1\right)}.
\end{align}
Similarly to the previous examples, this family contains DPPs covering a wide range of repulsiveness measured by $\eta_n$, and as $n \rightarrow \infty$, they are more repulsive in the sense that $\mathds{E}[\eta_n(\R^n)]$ decays slower. Indeed, 
 \begin{align*}
\mathds{E}[\eta_n(\R^n)] &= \frac{1}{e^{n\rho}} \int_{\R^n} K\added{_n}(x)^2 dx   = \added{e^{n\rho}} \frac{(2\pi)^{n/2}\alpha^{n}}{(\sigma + n)^{n/2} \Gamma\left(\frac{n}{2}\right)}  \frac{\Gamma\left( \frac{\sigma + n + 2}{2}\right)^2\Gamma\left( \frac{n}{2}\right) \Gamma(\sigma +1)}{ \Gamma\left(\frac{\sigma}{2} + 1\right)^2 \Gamma\left(\sigma + \frac{n}{2} + 1 \right)} \\
&= e^{n\rho} \frac{(2\pi)^{n/2}\alpha^{n}}{(\sigma + n)^{n/2} }  \frac{\Gamma(\sigma +1)\Gamma\left( \frac{\sigma}{2} + \frac{n}{2} + 1 \right)^2}{ \Gamma\left(\frac{\sigma}{2} + 1\right)^2 \Gamma\left(\sigma + \frac{n}{2} + 1 \right)},
\end{align*}
and by the upper bound \eqref{a_bnd_bessel},
\begin{align*}
\mathds{E}[\eta\added{_n}(\R^n)]  &<  \frac{\Gamma(\sigma +1)\Gamma\left( \frac{\sigma}{2} + \frac{n}{2} + 1 \right)}{ \Gamma\left(\frac{\sigma}{2} + 1\right) \Gamma\left(\sigma + \frac{n}{2} + 1 \right)}.
\end{align*}
By Stirling's formula, as $n \rightarrow \infty$, $\frac{\Gamma(\sigma +1)\Gamma\left( \frac{\sigma}{2} + \frac{n}{2} + 1 \right)}{ \Gamma\left(\frac{\sigma}{2} + 1\right) \Gamma\left(\sigma + \frac{n}{2} + 1 \right)} = O(n^{-\sigma/2})$.

These DPPs do not satisfy the conditions of Proposition \ref{thresh}, and so the concentration of the first moment measure does not occur, contrary to the first two families presented. However, the repulsive measure does not reach past the $\sqrt{n}$ scale in the sense of the following proposition.

\begin{prop}\label{bessel} Let $\rho \in \R$, $\alpha > 0$, and \added{$\sigma > 0$.} For each $n$, let $\Phi_n \sim DPP(K_n)$ in $\R^n$ \added{with $K_n$ given by \eqref{e:Bessel_K}.} 
Then, for any $\beta > \frac{1}{2}$ and $R > 0$,
\begin{align*}
\lim_{n \to \infty} \frac{\mathds{E}[\eta_n(\R^n \backslash B_n(Rn^{\beta}))]}{\mathds{E}[\eta_n(\R^n)]} = 0.
\end{align*}
\end{prop}

\subsection{Normal Variance Mixture Models}

Another class of DPPs described in \cite{Moller} are those with normal-variance mixture kernels. Let $\Phi_n \sim DPP(K_n)$ be a normal-variance mixture DPP in $\R^n$ with intensity $e^{n\rho}$ for $\rho \in \R$, i.e., let
\begin{align*}
K_n(x) = e^{n\rho} \frac{\mathds{E}[W^{-n/2}e^{-|x|^2/(2W)}]}{\mathds{E}[W^{-n/2}]},\qquad x \in \R^n,
\end{align*}
for some non-negative real-valued random variable $W$ such that $\mathds{E}[W^{-n/2}] < \infty$.
From \cite{Moller}, the bound $0 \leq \hat{K} <1$ translates to the following bound on the intensity:
\begin{align}\label{e:nv_int_bnd}
e^{n\rho} < \mathds{E}[W^{-n/2}]/(2\pi)^{n/2}.
\end{align}
If $\sqrt{2W} = \alpha$, this is known as the Gaussian DPP model. If $W \sim \text{Gamma}(\nu + \frac{n}{2}, 2\alpha^2)$, this is called the Whittle-Mat\'{e}rn model. The Cauchy model is given when $\frac{1}{W} \sim \text{Gamma}(\nu, 2\alpha^{-2})$. In both cases $\nu > 0$ and $\alpha > 0$ are parameters.

This family of DPPs does not cover a wide range of repulsiveness like the previous families. Indeed, for any random variable $W$ in $\R^+$ such that $\mathds{E}[W^{-\frac{n}{2}}] < \infty$, Parseval's theorem, Jensen's inequality, \eqref{e:nv_int_bnd}, and Fubini's theorem imply
\begin{align*}
\mathds{E}[\eta_n(\R^n)] &= \frac{1}{e^{n\rho}} \int_{\R^n} \hat{K}_n(x)^2 {\mathrm d}x 
= \frac{1}{e^{n\rho}} \int_{\R^n} \left(e^{n\rho} \tfrac{(2\pi)^{\frac{n}{2}}}{\mathds{E}\left[W^{-\frac{n}{2}}\right]} \mathds{E}\left[e^{-2\pi^2 |x|^2 W}\right]\right)^2 {\mathrm d}x \\
&\leq \tfrac{(2\pi)^{\frac{n}{2}}}{\mathds{E}\left[W^{-\frac{n}{2}}\right]} \int_{\R^n} \mathds{E}\left[e^{-4\pi^2 |x|^2 W}\right] {\mathrm d}x \\
&=   \tfrac{(2\pi)^{\frac{n}{2}}}{\mathds{E}[W^{-\frac{n}{2}}]} \mathds{E}\left( (4\pi W)^{-\frac{n}{2}} \mathds{E}\left[ (4\pi W)^{\frac{n}{2}} \int_{\R^n} e^{-4\pi^2 |x|^2 W} {\mathrm d}x \,\bigg| \, W \right] \right) = 2^{-\frac{n}{2}}.
\end{align*}

Is it difficult to make further general statements about this class because the behavior of the sequence $\frac{|X_n|}{\sqrt{n}}$ depends greatly on the distribution of the $\R^{+}$-valued random variable $W$. The rest of the section will describe results for specific models in this class.

Consider a sequence of normal-variance mixture DPPs all associated with the same random variable $W$. If $W$ is a constant $\alpha$, the random variables $X_n$ become multivariate Gaussian vectors with mean zero and variance depending on $\alpha$. The scaled norms of these vectors are well-known to satisfy a LDP since the coordinates are independent. This also corresponds to a Laguerre-Gaussian DPP with parameter $m = 2$. 

There is also a subclass of the normal-variance mixture models that satisfy Proposition \ref{logconcave}. In \cite{Yu}, it is proved that if $W$ has a log-concave density, then the normal-variance mixture distribution is log-concave. This implies that $K_n^2$ is log-concave, and thus if condition (\ref{e:meanconv}) holds, the conclusion of Proposition \ref{logconcave} holds. 
\added{Since the Gamma distribution for shape parameter $\nu$ greater than 1 is log-concave and $\nu + \frac{n}{2} \geq 1$ for large $n$, Whittle-Mat\'{e}rn DPPs are an example from this subclass and exhibit an $R^*$ as shown in the following proposition.

\begin{prop}\label{Whittle}
For each $n$, let $\Phi_n \sim DPP(K_n)$ be a Whittle-Mat\'{e}rn model in $\R^n$ with intensity $e^{n\rho}$ and parameters $\nu > 0$ and $\alpha > 0$, i.e., let 
\begin{align}\label{e:whittle_K}
K_n(x) = e^{n\rho}\frac{2^{1 - \nu}}{\Gamma(\nu)}\frac{|x|^{\nu}}{\alpha^{\nu}}\mathbb{K}_{\nu}\left(\frac{|x|}{\alpha}\right),\qquad x \in \R^n,
\end{align}
where $\alpha < \frac{\Gamma(\nu)^{\frac{1}{n}}}{\Gamma\left(\nu + \frac{n}{2}\right)^{\frac{1}{n}}2\sqrt{\pi} e^{\rho}}$ and $\mathbb{K}_{\nu}$ is the modified Bessel kernel of the second kind. Then, for $R^* := \frac{\alpha}{2}$,
\begin{align*}
\lim_{n \rightarrow \infty} \frac{\mathds{E}[\eta_n(B_n(\sqrt{n}R))]}{\mathds{E}[\eta_n(\R^n)]} = \begin{cases} 0, & R < R^* \\ 1, & R > R^*. \end{cases}
\end{align*}
\end{prop}}

\added{\begin{rem}
The upper bound on $\alpha$ needed for existence implies that for all $\nu$,
\begin{align*}
R^* := \frac{\alpha}{2} < \frac{\Gamma(\nu)^{\frac{1}{n}}}{\Gamma\left(\nu + \frac{n}{2}\right)^{\frac{1}{n}}4\sqrt{\pi} e^{\rho}} < \frac{1}{\sqrt{2\pi e} e^{\rho}} := \tilde{R},
\end{align*}
since $\left(\frac{\Gamma(\nu)}{\Gamma\left(\nu + \frac{n}{2}\right)}\right)^{\frac{1}{n}} \leq 1$ and $4 > \sqrt{2e}$. Thus, for these models, $R^*$ never reaches past the nearest neighbor threshold $\tilde{R}$.
\end{rem}}

Finally, the following proposition shows that the Cauchy models satisfy the conditions of Proposition \ref{thresh} if the $\alpha$ parameter grows appropriately with $n$. 

\begin{prop}\label{Cauchy} For each $n$, let $\Phi_n \sim DPP(K_n)$ be a Cauchy model in $\R^n$ with intensity $e^{n\rho}$ and parameters $\nu > 0$ and $\alpha_n > 0$, i.e., let 
\begin{align*}
K_n(x) = \frac{e^{n \rho}}{(1 + |\frac{x}{\alpha_n}|^2)^{\nu + \frac{n}{2}}}, \qquad x \in \R^n.
\end{align*}
Assume $\alpha_n \sim \alpha n^{1/2}$ as $n \rightarrow \infty$ for some $\alpha > 0$ such that $\alpha_n <  \frac{\Gamma(\nu + \frac{n}{2})^{\frac{1}{n}}}{\sqrt{\pi}e^{\rho}\Gamma(\nu)^{\frac{1}{n}}}$ for each $n$. Then, for $R^* := \alpha$,
\begin{align*}
\lim_{n \rightarrow \infty} \frac{\mathds{E}[\eta_n(B_n(\sqrt{n}R))]}{\mathds{E}[\eta_n(\R^n)]} = \begin{cases} 0, & R < R^* \\ 1, & R > R^*. \end{cases}
\end{align*}
\end{prop}

\begin{rem}
The upper bound on $\alpha_n$ has the following asymptotic expansion as $n \rightarrow \infty$:
\begin{align*}
\alpha_n < \frac{\Gamma(\nu + \frac{n}{2})^{\frac{1}{n}}}{\sqrt{\pi}e^{\rho}\Gamma(\nu)^{\frac{1}{n}}} \sim \frac{n^{1/2}}{\sqrt{2\pi e} e^{\rho}}.
\end{align*}
Thus, if $\alpha_n \sim \alpha n^{\frac{1}{2}}$, the reach of repulsion has the upper bound
\begin{align*}
 R^* := \alpha < \frac{1}{\sqrt{2\pi e} e^{\rho}}.
 \end{align*}
 This upper bound is precisely the threshold $\tilde{R}$ for the nearest neighbor function, and so unlike in the case of Laguerre-Gaussian DPPs and power exponential DPPs, the reach of repulsion $R^*$ for a sequence of Cauchy models with fixed parameter $\nu$ cannot reach past $\tilde{R}$.
\end{rem}

\section{Application to determinantal Boolean models in the Shannon regime}\label{applications}
Poisson Boolean models in the Shannon regime were studied in \cite{Venkat}, and the degree threshold results can be extended to Laguerre-Gaussian DPPs using Proposition \ref{Laguerre_rate}.
 
The setting is the following: Consider a sequence of stationary DPPs $\Phi_n$, indexed by dimension, where $\Phi_n \sim DPP(K_n)$ in $\mathds{R}^n$. Assume that for each $n$, $K_n$ is continuous, symmetric, and $0 \leq \hat{K}_n < 1$. Let the intensity of $\Phi_n$ be $K_n(0) = e^{n\rho}$.
Let $\Phi_n = \sum_k \delta_{T_n^{(k)}}$ and $R > 0$. 
Then, consider the sequence of particle processes \cite{Weil}, called determinantal Boolean models,
\begin{align*}
\mathcal{C}_n = \bigcup_k B\added{_n}\left(T_n^{(k)},\frac{ \sqrt{n}R}{2}\right).
\end{align*} 
The degree of each model is the expected number of balls that intersect the ball centered at zero under the reduced Palm distribution, i.e., $\mathds{E}[\Phi_n^{0,!}(B_n(\sqrt{n}R))]$. In the case when $\Phi_n$ is Poisson,
$\mathds{E}^{0,!}[\Phi_n(B(\sqrt{n}R))] = \mathds{E}[\Phi_n(B(\sqrt{n}R))]$ by Slivnyak's theorem, and
\begin{align*}
\lim_{n \rightarrow \infty} \frac{1}{n} \ln \mathds{E}^{0,!}[\Phi_n(B_n(\sqrt{n}R))] = \rho + \frac{1}{2} \log 2 \pi e + \log R.
\end{align*} 
To extend this result to DPPs, it is needed that as $n \rightarrow \infty$,
\begin{center}
$\mathds{E}[\Phi_n^{0,!}(B_n(\sqrt{n}R))] \sim  \mathds{E}[\Phi_n(B_n(\sqrt{n}R))]$. 
\end{center}
\added{Note that this would be impossible for a repulsive point process like the Mat\'ern hardcore process, since $\mathds{E}[\Phi_n^{0,!}(B_n(R_n))] = 0$ for all $R_n $ less than the hardcore radius.} 

However, for DPPs, notice that
\begin{align*}
\frac{\mathds{E}[\Phi_n^{0,!}(B_n(\sqrt{n}R))]}{\mathds{E}[\Phi_n(B_n(\sqrt{n}R))] } = 1 -\frac{\mathds{E}[\eta_n(B_n(\sqrt{n}R))]}{\mathds{E}[\Phi_n(B_n(\sqrt{n}R))] }.
\end{align*}
Thus, if $\frac{\mathds{E}[\eta_n(B_n(\sqrt{n}R))]}{\mathds{E}[\Phi_n(B_n(\sqrt{n}R))] } \rightarrow 0$ as $n \rightarrow \infty$, then the degree of the determinantal Boolean model has the same asymptotic behavior as the Poisson Boolean model. 

In the case of Laguerre-Gaussian kernels, this is the case, and the earlier results even provide the rate at which the quantity goes to zero, which exhibits a threshold at $R^*$ as is expected. 

\begin{prop}\label{gauss_beta2} Let $m \in \mathds{N}$ and $\rho \in \R$. For each $n$, let $\Phi_n \sim DPP(K_n)$ in $\R^n$ where 
\begin{align*}
K_n(x) = \frac{e^{n\rho}}{\left(\begin{smallmatrix} m - 1 + n/2 \\ m - 1 \end{smallmatrix}\right)} L_{m-1}^{n/2}\left(\frac{1}{m} \left| \frac{x}{\alpha} \right|^2\right) e^{-\frac{|x/\alpha|^2}{m}},
\end{align*}
and $\alpha$ is a parameter such that $0 < \alpha < \frac{1}{\sqrt{m\pi} e^{\rho}}$.  Then,
\begin{align*}
\lim_{n \rightarrow \infty} - \frac{1}{n} \ln \frac{\mathds{E}[\eta_n(B_n(\sqrt{n}R))]}{\mathds{E}[\Phi_n(B_n(\sqrt{n}R))] }  = \begin{cases}  \frac{2R^2}{\alpha^2m}, & 0 < R < \sqrt{m}\frac{\alpha}{2} \\  \frac{1}{2} +\log 2 - \log \alpha - \frac{1}{2} \log m + \log R, & R > \sqrt{m}\frac{\alpha}{2}. 
\end{cases}
\end{align*}
\end{prop}

\section{Conclusion}\label{conclusion}

By examining a measure of repulsiveness of DPPs, this paper provides insight into the high dimensional behavior of different families of DPP models. \added{Most of the families of DPPs presented in this paper have a global measure of repulsion decreasing to zero as dimension increases, indicating that they become more and more similar to Poisson point processes in high dimensions by \eqref{e:coupling_ineq}.} However, the reach of the small repulsive effect can still be quantified. By making a connection between the kernel of the DPP and the concentration in high dimensions of the norm of a random vector, we have shown under certain conditions that there exists a distance on the $\sqrt{n}$ scale at which the repulsive effect of a point of the DPP model is strongest as $n \rightarrow \infty$. It has been illustrated that some families of DPPs exhibit this reach of repulsion and some do not. The results are summarized in Table \ref{ex_summary}.

Many questions remain concerning the range of possible repulsive behavior of DPPs in high dimensions. First, the results can be extended to scalings other than the Shannon regime in the following way. Assumption \added{\eqref{conv}} in Proposition \ref{thresh} can be generalized to the assumption that for some sequence $b_n$, \added{$\frac{|X_n|}{b_n} \rightarrow R^*$} as $n \rightarrow \infty$. If $b_n \neq O(n\added{^\frac{1}{2}})$, the result holds for a different scaling than the Shannon regime, and the repulsiveness is strongest near $R^*b_n$ in high dimensions. While this is precisely what is shown not to happen for the Bessel-type DPPs \added{if $\sigma > 0$}, \added{examples of this generalization for $b_n = o(n)$ can be obtained from the power exponential DPPs when $\alpha_n = o(n^{\frac{1}{\nu} - \frac{1}{2}})$}. However, as noted in the introduction, any distance scaling smaller than $\sqrt{n}$ will not reach the regime where the expected number of points goes to infinity as dimension grows. Thus, this scaling appears less interesting. It would be interesting to find a family of DPPs that exhibits the concentration for $b_n \gg \sqrt{n}$. 

For all of the DPPs studied in this paper, $\mathds{E}[\eta_n(\R^n)] \rightarrow 0$ as $n \rightarrow \infty$. This is not always the case. For instance, there exists a class of DPPs such that  for $c \in (0,1)$, $\mathds{E}[\eta_n(\R^n)] = c$ for all $n$. Indeed, let $K_n \in L^2(\R^n)$ be such that its Fourier transform is
\begin{equation}\label{mostrepulsive}
\hat{K}_n(\xi) = \sqrt{c}1_{B_n(r_{n})}(\xi), \xi \in \R^n,
\end{equation}
where $r_{n} \in \R^{+}$ is such that $\text{Vol}(B_n(r_{n})) = K_n(0)$. Then,
\begin{equation*}
\mathds{E}[\eta_n(\R^n)] = \frac{1}{K_n(0)} \int_{\R^n} K_n(x)^2 {\mathrm d}x =  \frac{1}{K_n(0)} \int_{\R^n} \hat{K}_n(\xi)^2 {\mathrm d}\xi =   \frac{c}{K_n(0)} \text{Vol}(B_n(r_{n}))  = c.
\end{equation*}
It would be useful to find a necessary and sufficient condition for $\mathds{E}[\eta_n(\R^n)]$ to converge to zero. However, if $\mathds{E}[\eta_n(\R^n)]$ does not converge to zero, this does not necessarily prevent $\mathds{P}(\eta_n(\R^n) = 0)$ from approaching one as $n$ goes to infinity. It would be interesting to find a class of DPPs where $\mathds{P}(\eta_n(\R^n) = 0)$ approaches some $c < 1$. 

There is an important class of stationary and isotropic DPPs that should be mentioned. Recall that \added{to ensure $\eta$ is well-defined, it is assumed that} the kernel $K$ associated with $\Phi$ satisfies $0 \leq \hat{K} < 1$. However, $\Phi$ still exists when $\hat{K}$ is allowed to attain the maximum value of one. For the models studied in this paper, it is the case when the parameter achieves its upper bound. In this case, we can still define the measure of repulsiveness (\ref{gamma}) even though it \added{may not be interpretable} as the intensity measure of a point process $\eta$. Replacing $\mathds{E}[\eta(B)]$ with $\int_B (1- g(x)){\mathrm d}x$ for $B \in \mathcal{B}(\R^n)$, the main results (Propositions \ref{thresh}, \ref{logconcave}, and \ref{LDP}) can be restated with the condition that $0 \leq \hat{K} \leq 1$.
In this case, the reach of repulsion $R^*$ is interpreted as the distance on the $\sqrt{n}$ scale at which the measure of repulsion is strongest.

A particularly interesting subclass of the DPPs described in the previous paragraph are the most repulsive stationary DPPs, introduced in the on-line supplementary material to \cite{Moller} (see \cite{Moller_Arxiv}). These DPPs maximize the measure of repulsiveness $\gamma$, and have a kernel $K$ such that $\hat{K}$ is defined as in (\ref{mostrepulsive}) but with $c = 1$.  For the most repulsive DPPs, $\gamma = 1$ in any dimension. In addition, for a sequence of DPPs $\{\Phi_n\}_{n \in \mathds{N}}$ where $\Phi_n$ is the most repulsive DPP in $\R^n$ with intensity $e^{n\rho}$, $X_n$ as defined in Proposition \ref{thresh} satisfies 
\begin{align*}
\mathds{E}[|X_n|^2] &= \int_{\R^n} |x|^2 \frac{K_n(x)^2}{||K_n||_2^2} {\mathrm d}x =  \frac{\Gamma\left(\frac{n}{2} + 1 \right)}{\pi^{n/2}} \int_{\R^n} |x|^2 \frac{J_{n/2}^2\left(2\sqrt{\pi} \Gamma \left(\frac{n}{2} + 1 \right)^{1/n} e^{\rho} |x| \right)}{|x|^n} {\mathrm d}x \\
&= n \int_0^{\infty} r J_{n/2}^2\left(2\sqrt{\pi} \Gamma \left(\frac{n}{2} + 1 \right)^{1/n} e^{\rho} r \right) {\mathrm d}r,
\end{align*}
where $J_{\nu}$ is the Bessel function of the first kind of order $\nu$ (see \cite{Lav}). By \cite[Eq.~1.17.13]{DLMF}, this integral does not converge, i.e., $|X_n|$ does not have a finite second moment.

\added{\begin{table}[h]
\centering
 \caption{Summary of Results}
  \begin{tabular}{ l | c | c | c | c}
    \hline
    DPP Class & $\mathds{E}[\eta_n(\R^n)]$ & $R^*$ & Rate type & $R^* > \tilde{R}$ \\ \hline
    Laguerre-Gaussian & $< 2^{-\frac{n}{2}}O(n^{m-1})$ & $\sqrt{m} \frac{\alpha}{2}$ & LDP & $\left(\frac{2}{e}\right)^{\frac{1}{2}} < e^{\rho} \sqrt{m \pi} \alpha < 1$ \\ 
    Power-Exponential  & $< 2^{-\frac{n}{\nu}}$ & $\alpha\frac{(2\nu)^{\frac{1}{\nu}}}{4\pi}$ & Chebychev & $\frac{2}{2^{\frac{1}{\nu}}e} < \frac{e^{\rho}\nu^{\frac{1}{\nu}}}{\sqrt{2\pi e}} \alpha < \frac{1}{e^{\frac{1}{\nu}}}$ \\ 
    Bessel-type & $ < O(n^{-\sigma/2})$ & N/A & N/A & N/A \\ 
    Whittle-Mat\'{e}rn & $< 2^{-\frac{n}{2}}$ & $\frac{\alpha}{2}$ & Log-concave & N/A \\ 
    Cauchy & $< 2^{-\frac{n}{2}}$ & $\alpha$ & Chebychev & N/A \\
    \hline
  \end{tabular}
  \label{ex_summary}
\end{table}}







\appendix

\section{Proof of (\ref{near_thresh})}\label{near_thresh_app}

For each $n$, let $\Phi_n \sim DPP(K_n)$ in $\R^n$ be stationary with intensity $K_n(0) = e^{n\rho}$. \added{From (\ref{Shannon}), there exists $\tilde{R} := \frac{1}{\sqrt{2\pi e} e^{\rho}}$ such that} 
\begin{align*}
\lim_{n \to \infty} \mathds{E}[\Phi_n(B_n(\sqrt{n}R))] = \begin{cases} 0, & R < \tilde{R} \\ \infty, & R > \tilde{R} .\end{cases}
\end{align*}
By Theorem \ref{palm},
\begin{align*}
\mathds{E}[\Phi_n(B_n(\sqrt{n}R))] - \mathds{E}[\Phi_n^{0,!}(B_n(\sqrt{n}R))] &= \frac{1}{e^{n\rho} }\int_{B_n(\sqrt{n}R)} K_n(x)^2 {\mathrm d}x
\end{align*}
Then, by Parseval's theorem and Theorem \ref{exist},
\begin{align*}
\frac{1}{e^{n\rho} }\int_{B_n(\sqrt{n}R)} K_n(x)^2 {\mathrm d}x &\leq \frac{1}{e^{n\rho} }\int_{\R^n} \hat{K}_n(\xi)^2 {\mathrm d}\xi \leq \frac{1}{e^{n\rho} }\int_{\R^n} \hat{K}_n(\xi) {\mathrm d}\xi = 1.
\end{align*}
Also, since $\frac{1}{e^{n\rho} }\int_{B_n(\sqrt{n}R)} K_n(x)^2 {\mathrm d}x \geq 0$, the following bounds hold:
\begin{align*}
\mathds{E}[\Phi_n(B_n(\sqrt{n}R))] - 1 \leq \mathds{E}[\Phi_n^{0,!}(B_n(\sqrt{n}R))] \leq \mathds{E}[\Phi_n(B_n(\sqrt{n}R))].
\end{align*}
Thus, the threshold remains the same for the reduced Palm expectation:
\begin{align*}
\lim_{n \to \infty} \mathds{E}[\Phi^{0,!}_n(B_n(\sqrt{n}R))] = \begin{cases} 0, & R < \tilde{R} \\ \infty, & R > \tilde{R} .\end{cases}\end{align*}
By the first moment inequality and Proposition 5.1 in \cite{Blasz}, one has the following bounds:
\begin{align*}
1 - \mathds{E}[\Phi^{0,!}_n(B(\sqrt{n}R))]  \leq \mathds{P}(\Phi^{0,!}_n(B_n(\sqrt{n}R)) = 0)  \leq \exp\left(-\mathds{E}[\Phi^{0,!}_n(B(\sqrt{n}R))]  \right).
\end{align*}
Thus, 
$\lim_{n \rightarrow \infty} \mathds{P}(\Phi^{0,!}_n(B_n(\sqrt{n}R)) > 0) = \begin{cases} 0, & R < \tilde{R} \\ 1 & R > \tilde{R}.\end{cases}$


\added{\section{Proof of Lemma \ref{firstmoment}}


By Theorem \ref{palm}, for any $B \in \mathcal{B}(\R^n)$,
\begin{align}\label{eta_density}
\mathds{E}[\eta_n(B)] &= \mathds{E}[\Phi_n(B)] - \mathds{E}[\Phi_n^{0,!}(B)] 
= \frac{1}{K_n(0)} \int_B K_n(x)^2 {\mathrm d}x,
\end{align}
i.e. the first moment measure of $\eta_n$ has a density with respect to Lebesgue measure equal to $\frac{1}{K_n(0)}K_n(x)^2$. Then by the monotone convergence theorem,
\begin{align*}
\mathds{E}[\eta_n(\R^n)] = \lim_{R \rightarrow \infty} \mathds{E}[\eta_n(B_n(R))] = \frac{1}{K_n(0)} \int_{\R^n} K_n(x)^2 {\mathrm d}x = \frac{||K_n||_2^2}{K_n(0)}.
\end{align*}
\added{Thus,} for all $B \in \mathcal{B}(\R^n)$,
\added{\begin{align*}
 \mathds{P}(X_n \in B) =  \int_B \frac{K_n(x)^2}{||K_n||_2^2} {\mathrm d}x = \frac{\mathds{E}[\eta_n(B)] }{ \mathds{E}[\eta_n(\R^n)]} .
 \end{align*}}}

\section{Proof of Main Results}
\subsection{Proof of Proposition \ref{thresh}}\label{thresh_app}

The assumption $\frac{|X_n|}{\sqrt{n}} \rightarrow R^*$ in probability means that for all $\varepsilon > 0$,
\[
\mathds{P}\left( \bigg| \frac{|X_n|}{\sqrt{n}} - R^* \bigg| > \varepsilon\right) \rightarrow 0, \text{ as } n \rightarrow \infty. 
\]
Now, assume $R  <  R^*$. Then, there exists $\varepsilon > 0$ such that $R = R^* -  \varepsilon $. Thus, 
\begin{align*}
\mathds{P}(|X_n| \leq \sqrt{n}R) 
= \mathds{P}\left(\frac{|X_n|}{\sqrt{n}} \leq R^* - \varepsilon \right) 
\leq \mathds{P}\left( \bigg|\frac{|X_n|}{\sqrt{n}} - R^* \bigg| > \varepsilon \right) \rightarrow 0 \text{ as } n \rightarrow \infty. 
\end{align*}
Second, assume $R > R^*$. Then, there exists $\varepsilon >  0$ such that  $R = R^* + \varepsilon$, and 
\begin{align*}
\mathds{P}(|X_n| \leq \sqrt{n}R) =  1 - \mathds{P}\left(\frac{|X_n|}{\sqrt{n}} > R^* + \varepsilon \right)  \geq 1 - \mathds{P}\left(\bigg|\frac{|X_n|}{\sqrt{n}} - R^* \bigg|>  \varepsilon \right) \rightarrow 1.
\end{align*}
Then, by Lemma \ref{firstmoment}, as $n \rightarrow \infty$,
\begin{align*}
\frac{\mathds{E}[\eta_n(B_n(\sqrt{n}R))]}{\mathds{E}[\eta_n(\R^n)]} = \mathds{P}\left(|X_n| \leq \sqrt{n}R \right) \rightarrow \begin{cases} 0, & R < R^* \\ 1, & R > R^*. \end{cases}
\end{align*}

\subsection{Proof of Proposition \ref{logconcave}}

\added{Since for all $n$, $\Phi_n$ is isotropic, $X_n$ as defined in Proposition \ref{thresh} has a radially symmetric density. Thus, $X_n$ has the same distribution as the product $R_n U_n$, where $R_n$ is equal in distribution to  $|X_n|$, $U_n$ is uniformly distributed on $\mathds{S}^{n-1}$, and $R_n$ and $U_n$ are independent. Letting $\sigma_n^2 = \mathds{E}|X_n|^2$ for each $n$, $\frac{\sqrt{n}}{\sigma_n}X_n$ then satisfies the conditions of Theorem \ref{thinshellthm} for each $n$. }Then, by Theorem \ref{thinshellthm}, for any $\delta > 0$, there exist absolute constants $C, c > 0$ such that
\begin{align*}
 \mathds{P}\left( \bigg|\frac{|X_n|}{\sigma_n} - 1 \bigg| \geq \delta \right) \leq C e^{-c \sqrt{n} \min(\delta, \delta^3) }.
\end{align*}
Now, let $\delta \in (0,1)$. By Lemma \ref{firstmoment}, 
\begin{align*}
\frac{\mathds{E}[\eta_n(B_n(\sigma_n(1 - \delta)))]}{\mathds{E}[\eta_n(\R^n)]} = \mathds{P}\left(\frac{|X_n|}{\sigma_n} \leq 1 - \delta \right) \leq C e^{-c \sqrt{n} \delta^3 },
\end{align*}
since $\min(\delta^3, \delta) = \delta^3$ for $\delta \in (0,1)$. Similarly, for any $\delta > 0$, 
\begin{align*}
 \frac{\mathds{E}[\eta_n(\R^n \setminus B_n(\sigma_n(1 + \delta)))]}{\mathds{E}[\eta_n(\R^n)]} = \mathds{P}\left(\frac{|X_n|}{\sigma_n} \geq 1 + \delta \right) \leq C e^{-c \sqrt{n}\min(\delta^3, \delta)}.
\end{align*} 
Now, assume $\frac{\sigma_n}{\sqrt{n}} \rightarrow R^* \in (0, \infty)$ as $n \rightarrow \infty$. For $R < R^*$, there exists $\varepsilon \in (0,1)$ such that $R = R^*(1 - \varepsilon)$. Then, for all $n$ large enough, $\frac{\sqrt{n}R^*}{\sigma_n} < \frac{1 - \frac{\varepsilon}{2}}{1 - \varepsilon}$ and 
\begin{align*}
\frac{\mathds{E}[\eta_n(B_n(\sqrt{n}R))]}{\mathds{E}[\eta_n(\R^n)]} &= \mathds{P}\left(|X_n| \leq \sqrt{n}R \right) = \mathds{P}\left(\frac{|X_n|}{\sigma_n} \leq \frac{\sqrt{n}R}{\sigma_n} \right) \\
&= \mathds{P}\left( \frac{|X_n|}{\sigma_n} \leq \frac{\sqrt{n}R^*(1 - \varepsilon)}{\sigma_n}\right) \leq  \mathds{P}\left( \frac{|X_n|}{\sigma_n} \leq 1 - \frac{\varepsilon}{2} \right)\\
&\leq \mathds{P}\left( \bigg|\frac{|X_n|}{\sigma_n} - 1\bigg| \geq \frac{\varepsilon}{2} \right) \leq C e^{-c \sqrt{n} ({\varepsilon}/{2})^3 }.
\end{align*}
Thus for all $R < R^*$, there exists a constant $C(\varepsilon(R)) = c {\varepsilon}^3/{2}^3$ such that
\begin{equation*}
\liminf_{n \rightarrow \infty} - \frac{1}{\sqrt{n}} \ln \frac{\mathds{E}[\eta_n(B_n(\sqrt{n}R))]}{\mathds{E}[\eta_n(\R^n)]} \geq C(\varepsilon(R)).
\end{equation*}
A similar argument gives that for all $R > R^*$, there exists $C(\varepsilon(R))$ such that 
\begin{equation*}
\liminf_{n \rightarrow \infty} - \frac{1}{\sqrt{n}} \ln \frac{\mathds{E}[\eta_n(\R^n \setminus B_n(\sqrt{n}R))]}{\mathds{E}[\eta_n(\R^n)]} \geq C(\varepsilon(R)).
\end{equation*}
This implies the threshold \eqref{threshold}.

\subsection{Proof of Proposition \ref{LDP}}

If $\frac{|X_n|}{\sqrt{n}}$ satisfies a large deviations principle with convex rate function $I$, then by definition,
\begin{align*}
- \inf_{r < R} I(r) \leq \liminf_{n \rightarrow \infty} \frac{1}{n} \ln \mathds{P}\left(\frac{|X_n|}{\sqrt{n}} \leq R \right)\leq \limsup_{n \rightarrow \infty} \frac{1}{n} \ln \mathds{P}\left(\frac{|X_n|}{\sqrt{n}} \leq R \right) \leq - \inf_{r \leq R} I(r). 
\end{align*} 
Thus, by Lemma \ref{firstmoment},
\begin{align*}
- \inf_{r < R} I(r) \leq \liminf_{n \rightarrow \infty} \frac{1}{n} \ln \frac{\mathds{E}[\eta_n(B_n(\sqrt{n}R))]}{\mathds{E}[\eta_n(\R^n)]}\leq \limsup_{n \rightarrow \infty} \frac{1}{n} \ln \frac{\mathds{E}[\eta_n(B_n(\sqrt{n}R))]}{\mathds{E}[\eta_n(\R^n)]} \leq - \inf_{r \leq R} I(r). 
\end{align*}
By the assumption that the rate function $I$ is strictly convex, there exists a unique $R^*$ such that $I(R^*) = 0$. Note that \added{$\inf_{\{r \leq R\}} I(r)$} is then zero for $R > R^*$. Thus,
\begin{align*}
\lim_{n \rightarrow \infty} \frac{\mathds{E}[\eta_n(B_n(\sqrt{n}R))]}{\mathds{E}[\eta_n(\R^n)]} =\begin{cases} 0, & R < R^* \\ 1, & R > R^*. \end{cases}
\end{align*}
Let $R < R^*$. If the rate function $I$ is continuous at $R$, then the above inequalities become equalities and 
\begin{align*}
\lim_{n \rightarrow \infty} - \frac{1}{n} \ln \frac{\mathds{E}[\eta_n(B_n(\sqrt{n}R))]}{\mathds{E}[\eta_n(\R^n)]} = I(R).
\end{align*}

\section{Proof of Lemma \ref{Laguerre}}

The proof shows that the sequence of random variables satisfies the conditions of the G\"{a}rtner-Ellis theorem (see \cite{LDP}). First,
\begin{equation*}
\mathds{E}[e^{s|X_n|^2}] = \frac{e^{2n\rho}}{\binom{m - 1 + n/2 }{m - 1 }^2 ||K_n||^2_2} \overbrace{\int_{\R^n} e^{-\left(\frac{2}{\alpha^2 m} - s \right)|x|^2} \left(L_{m-1}^{n/2} \left(\tfrac{1}{m} \left|\tfrac{x}{\alpha}\right|^2\right)\right)^2 {\mathrm d}x}^\text{I(s)}.
\end{equation*}
Writing out the polynomial, the integral $I$ above becomes
\begin{equation*}
I(s) = \sum_{k, j=0}^{m-1} \binom{ m - 1 + n/2 }{m - 1 - k } \binom{ m - 1 + n/2 }{m - 1 - j } \frac{(-1)^{k+j}}{k!j! (m\alpha^2)^{k+j}} \int_{\R^n} e^{-\left(\frac{2}{\alpha^2 m} - s \right)|x|^2}  |x|^{2k+2j} {\mathrm d}x.
\end{equation*}
A quick calculation shows that for $a > 0$,
\begin{equation}\label{important_integral}
\int_{\R^n} e^{-a|x|^2}|x|^{b} {\mathrm d}x = \frac{\pi^{n/2}}{a^{\frac{n + b}{2}}} \frac{\Gamma\left(\frac{n}{2} + \frac{b}{2}\right)}{\Gamma\left(\frac{n}{2}\right)}.
\end{equation}
 Then, if $s < \frac{2}{\alpha^2m}$,
\begin{equation*}
I(s) =  \frac{\pi^{n/2}}{\left(\frac{2}{\alpha^2 m } - s\right)^{\frac{n}{2}}\Gamma\left(\frac{n}{2}\right)} \sum_{k, j =0}^{m-1} \binom{ m - 1 + n/2}{ m - 1 - k}\binom{m - 1 + n/2}{ m - 1 - j} \frac{(-1)^{k+j}\Gamma\left(\frac{n}{2} + k + j \right)}{k!j!\left(2 - s m \alpha^2\right)^{k + j}},
\end{equation*}
and $I(s) = \infty$ otherwise.
For each $k, j \in \mathds{N}$,
\begin{align}\label{sim_lag}
&\binom{m - 1 + n/2}{m - 1 - k } \binom{ m - 1 + n/2}{ m - 1 - j } \Gamma\left(\frac{n}{2} + k + j \right)  \nonumber \\
&\qquad \sim \added{\frac{1}{(m-1-k)!(m-1-j)!}} \left(\frac{n}{2}\right)^{2m-2} \Gamma\left(\frac{n}{2}\right), 
\end{align}
as $n \rightarrow \infty$. So, $I(s)$ has the following asymptotic expansion for $s < \frac{2}{\alpha^2m}$ as $n \rightarrow \infty$:
\begin{align*}
I(s) \sim \frac{\pi^{n/2}}{\left(\frac{2}{\alpha^2 m } - s\right)^{\frac{n}{2}}}  \left(\frac{n}{2}\right)^{2m-2} \sum_{k=0}^{m-1} \added{ \sum_{j=0}^{m-1} } \frac{(-1)^{k+j}}{k!j!\added{(m-1-k)!(m-1-j)!}}   \frac{1}{\left(2 - s m \alpha^2\right)^{k + j}}.
\end{align*}
By (\ref{Lag_norm}) and (\ref{sim_lag}),
\begin{equation}\label{e:kasymp}
\frac{1}{e^{2n\rho}} ||K_n||_2^2  \sim \frac{\alpha^n}{ \binom{m - 1 + n/2}{m - 1}^2} \left(\frac{m\pi}{2}\right)^{\frac{n}{2}} \left(\frac{n}{2}\right)^{2m-2} \sum_{k, j=0}^{m-1} \frac{(-1)^{k+j}}{k!j!\added{(m-1-k)!(m-1-j)!}} \frac{1}{2^{k + j}},
\end{equation}
and hence,
\begin{align*}
\mathds{E}[e^{s|X_n|^2}]  &\sim \left(1 - \frac{s \alpha^2 m}{2}\right)^{-\frac{n}{2}}\left( \frac{\sum_{k, j=0}^{m-1}  \frac{(-1)^{k+j}}{k!j!\added{(m-1-k)!(m-1-j)!}}   \frac{1}{\left(2 - s m \alpha^2\right)^{k + j}} }{ \sum_{k, j=0}^{m-1}   \frac{(-1)^{k+j}}{k!j!\added{(m-1-k)!(m-1-j)!}}   \frac{1}{2^{k + j}}} \right),
\end{align*}
as $n \rightarrow \infty$. Thus,
\begin{align*}
\Lambda (s) = \lim_{n \rightarrow \infty} \frac{1}{n} \log \mathds{E}[e^{s|X_n|^2}]  = - \frac{1}{2} \log \left(1 - \frac{s \alpha^2 m}{2}\right) \text{ if } s < \frac{2}{\alpha^2m},
\end{align*}
and is infinite otherwise. It is clear that $0 \in (D(\Lambda))^{\circ}$, where $D(\Lambda) = \{ s \in \R : \Lambda(s) < \infty\}$. Thus, the G\"{a}rtner-Ellis conditions are satisfied. 
The rate function for the LDP is computed with the optimization
\begin{align*}
\Lambda^*(x) &= \sup_{\lambda \in \R} \left[ x \lambda - \Lambda(\lambda) \right] = \sup_{\lambda \in \R} \left[ x \lambda +  \frac{1}{2} \log \left(1 - \frac{\lambda \alpha^2 m}{2}\right) \right]. 
\end{align*}
Then, since
\begin{align*}
0 = \frac{d}{d\lambda} \left[ x \lambda +  \frac{1}{2} \log \left(1 - \frac{\lambda \alpha^2 m}{2}\right) \right] & = x - \frac{\alpha^2 m}{4 - 2\alpha^2 m \lambda} \text{ if and only if }
\lambda = \frac{2}{\alpha^2 m} - \frac{1}{2x},
\end{align*}
the rate function is
\begin{align*}
\Lambda^*(x) = x \left(\frac{2}{\alpha^2 m} - \frac{1}{2x}\right) +  \frac{1}{2} \log \left(1 - \frac{\left(\frac{2}{\alpha^2 m} - \frac{1}{2x}\right) \alpha^2 m}{2}\right) = \frac{2x}{\alpha^2 m} - \frac{1}{2} +  \frac{1}{2} \log \left(\frac{\alpha^2 m}{4x}\right).
\end{align*}
Then by the contraction principle (see \cite{LDP}), the sequence $\frac{|X_n|}{\sqrt{n}}$ satisfies an LDP with rate function
\begin{align*}
\Lambda^*(x) = \frac{2x^2}{\alpha^2 m} - \frac{1}{2} +  \frac{1}{2} \log \left(\frac{\alpha^2 m}{4x^2}\right).
\end{align*}
Note that $\Lambda^*(x) = 0$ if and only if $x = \sqrt{m}\frac{\alpha}{2}$\added{, implying $\frac{|X_n|}{\sqrt{n}} \rightarrow \sqrt{m} \frac{\alpha}{2}$ in probability.}

\section{Proof of Proposition \ref{Laguerre_rate}}

\begin{proof}
For each $n$, let $X_n$ be a random vector in $\R^n$ with density $\frac{K_n^2}{||K_n||_2^2}$. By Lemma \ref{Laguerre}, for $R < \sqrt{m}\frac{\alpha}{2}$,
\begin{align*}
 \lim_{n \rightarrow \infty} - \frac{1}{n} \log \mathds{P} \left( \frac{|X_n|}{\sqrt{n}} \leq R\right) = \frac{2R^2}{\alpha^2 m} - \frac{1}{2} +  \frac{1}{2} \log \left(\frac{\alpha^2 m}{4R^2}\right).
\end{align*}
Then by \eqref{e:kasymp}, as $n \rightarrow \infty$, 
\begin{equation}
\mathds{E}[\eta_n(\R^n)] = \frac{1}{e^{n\rho}} ||K_n||_2^2  \sim \left(\frac{e^{2\rho} \alpha^2 m\pi}{2}\right)^{\frac{n}{2}} \sum_{k, j=0}^{m-1} \frac{(-1)^{k+j}}{k!j!\added{(m-1-k)!(m-1-j)!}} \frac{1}{2^{k + j}},
\end{equation}
Thus, by Lemma \ref{firstmoment},
\begin{align*}
&\lim_{n \rightarrow \infty} - \frac{1}{n} \log  \mathds{E}[\eta_n(B_n(\sqrt{n}R))] =  \lim_{n \rightarrow \infty} - \frac{1}{n} \log \mathds{E}[\eta_n(\R^n)] + \lim_{n \rightarrow \infty} - \frac{1}{n} \log \mathds{P} \left( \frac{|X_n|}{\sqrt{n}} \leq R\right) \\
&= \begin{cases} - \rho - \log \alpha - \frac{1}{2} \log \left(\frac{m\pi}{2}\right) + \left(\frac{2R^2}{\alpha^2 m} - \frac{1}{2} +  \frac{1}{2} \log \left(\frac{\alpha^2 m}{4R^2}\right)\right), & 0 < R < \sqrt{m} \frac{\alpha}{2} \\
 - \rho - \log \alpha - \frac{1}{2} \log \left(\frac{m\pi}{2}\right),  & R >  \sqrt{m} \frac{\alpha}{2} \end{cases} \\
 &= \begin{cases} - \rho - \frac{1}{2} \log 2\pi e  + \frac{2R^2}{\alpha^2m} - \log R, & 0 < R < \sqrt{m}\frac{\alpha}{2} \\ - \rho - \log \alpha - \frac{1}{2} \log \frac{m\pi}{2} , & R > \sqrt{m}\frac{\alpha}{2}.
\end{cases}
\end{align*}
\end{proof}

\section{Proof of Lemma \ref{scale}}

Since for all $n$, $\hat{K}_n \in C^2(\R^n)$, Parseval's theorem implies
\begin{align}\label{Enorm}
\mathds{E}[|X_n|^2] & = \frac{1}{||K_n||_2^2} \int_{\R^n} |x|^2K_n(x)^2 {\mathrm d}x
 = \frac{1}{||\hat{K}_n||_2^2} \int_{\R^n} -\frac{\triangle \hat{K}_n(\xi)}{(2\pi)^2} \hat{K}_n(\xi){\mathrm d}\xi.
\end{align}
To compute the Laplacian of $\hat{K}$, we first see that for each $i$,
\begin{align*}
\frac{\partial^2}{\partial x_i^2} e^{- | \alpha x |^{\nu}} 
&=\frac{\partial}{\partial x_i}( - \nu \alpha^{\nu} x_i | x|^{\nu - 2} e^{- |\alpha x |^{\nu}}) \\
&= - \nu \alpha^{\nu} | x|^{\nu - 2} e^{- | \alpha x |^{\nu}} - \nu \alpha^{\nu} x_i \left( \tfrac{\partial}{\partial x_i} | x|^{\nu - 2} \right) e^{- | \alpha x |^{\nu}} + (\nu \alpha^{\nu} x_i | x|^{\nu - 2})^2 e^{- | \alpha x |^{\nu}} \\
&=  e^{- | \alpha x |^{\nu}} \left( - \nu \alpha^{\nu} | x|^{\nu - 2} - \nu(\nu -2) \alpha^{\nu} x_i^2 | x|^{\nu - 4} + \nu^2 \alpha^{2\nu} x_i^2 | x|^{2\nu - 4} \right) \\
&=  e^{- | \alpha x |^{\nu}} \left(x_i^2(\nu^2 \alpha^{2\nu}| x|^{2\nu - 4} - \nu(\nu -2) \alpha^{\nu} | x|^{\nu - 4}) - \nu \alpha^{\nu} | x|^{\nu - 2} \right).
\end{align*}
Then,
\begin{align*}
\triangle e^{- | \alpha x |^{\nu}}  &= \sum_{i=1}^n \frac{\partial^2}{\partial x_i^2} e^{- | \alpha x |^{\nu}}\\
&= \sum_{i=1}^n e^{- | \alpha x |^{\nu}} \left(x_i^2(\nu^2 \alpha^{2\nu}| x|^{2\nu - 4} - \nu(\nu -2) \alpha^{\nu} | x|^{\nu - 4}) - \nu \alpha^{\nu} | x|^{\nu - 2} \right) \\
&= e^{- | \alpha x |^{\nu}} \left(|x|^2(\nu^2 \alpha^{2\nu}| x|^{2\nu - 4} - \nu(\nu -2) \alpha^{\nu} | x|^{\nu - 4}) - n\nu \alpha^{\nu} | x|^{\nu - 2} \right) \\
&= e^{- | \alpha x |^{\nu}} \left(\nu^2 \alpha^{2\nu}| x|^{2\nu - 2} - ( \nu(\nu -2) \alpha^{\nu} + n\nu \alpha^{\nu} ) | x|^{\nu - 2}  \right).
\end{align*}
Thus by (\ref{Enorm}) and (\ref{powerexp_global}),
\begin{align*}
&\mathds{E}[|X_n|^2] = \frac{\Gamma(\frac{n}{2} + 1) \alpha_n^n2^{\frac{n}{\nu}}}{4\pi^2\pi^{n/2}\Gamma(\frac{n}{\nu} + 1)}  \int_{\R^n} e^{- 2| \alpha_n x |^{\nu}} \left( ( \nu(\nu -2) \alpha_n^{\nu} + n\nu \alpha_n^{\nu} ) | x|^{\nu - 2}  - \nu^2 \alpha_n^{2\nu}| x|^{2\nu - 2} \right) {\mathrm d}x \\
&\qquad = \frac{\Gamma(\frac{n}{2} + 1) \alpha_n^{n+\nu} 2^{\frac{n}{\nu}}\nu }{4\pi^2\pi^{\frac{n}{2}}\Gamma(\frac{n}{\nu} + 1)} \left[  (\nu -2 + n)\int_{\R^n} | x|^{\nu - 2} e^{- 2| \alpha_n x |^{\nu}}   {\mathrm d}x  - \nu \alpha_n^{\nu}\int_{\R^n} e^{- 2| \alpha_n x |^{\nu}} | x|^{2\nu - 2} dx \right].
\end{align*}
Then, using (\ref{important_integral}),
\begin{align*}
\mathds{E}[|X_n|^2] &= n\frac{\alpha_n^{n+\nu}2^{\frac{n}{\nu}}\nu }{4\pi^2\Gamma(\frac{n}{\nu} + 1)} \Bigg[ -\frac{ \nu \alpha_n^{\nu} \Gamma(\frac{n + 2\nu - 2}{\nu})}{\nu 2^{(n + 2\nu - 2)/\nu}\alpha_n^{n + 2\nu - 2}} + \frac{(\nu -2 + n)\Gamma(\frac{n + \nu - 2}{\nu})}{\nu 2^{(n + \nu - 2)/\nu}\alpha_n^{n + \nu - 2}} \Bigg] \\
&=  n  \frac{2^{2/\nu}\alpha_n^2}{4\pi^2\Gamma(\frac{n}{\nu} + 1)} \left[   \frac{(\nu -2 + n)}{2}\Gamma\left(\frac{n - 2}{\nu} + 1\right) - \frac{\nu}{4} \Gamma\left(\frac{n - 2}{\nu} + 2\right) \right] \\
&= n \frac{ 2^{2/\nu} \alpha_n^{2}\Gamma\left(\frac{n-2}{\nu} + 1\right)}{4\pi^2\Gamma(\frac{n}{\nu} + 1)} \left[ \frac{n}{4} + \frac{\nu}{4} - \frac{1}{2} \right].
\end{align*}
By the asymptotic formula for the Gamma function, as $n \rightarrow \infty$,
\begin{align*}
\mathds{E}[|X_n|^2]  &\sim n\frac{\alpha_n^2 2^{\frac{2}{\nu}}}{4\pi^2} \left(\sqrt{\frac{\nu}{2\pi n}} \left( \frac{\nu e}{n} \right)^{\frac{n}{\nu}} \right)\left(\sqrt{\frac{2\pi (n-2)}{\nu}} \left( \frac{n-2}{\nu e} \right)^{\frac{(n-2)}{\nu}} \right)\left[ \frac{n}{4} + \frac{\nu}{4} - \frac{1}{2} \right] \\
&= n \frac{\alpha_n^2 2^{2/\nu}}{4\pi^2}  \frac{\sqrt{n-2}}{\sqrt{n}} \left( 1- \frac{2}{n} \right)^{\frac{n}{\nu}} \left( \frac{n-2}{\nu e} \right)^{-\frac{2}{\nu}}\left[ \frac{n}{4} + \frac{\nu}{4} - \frac{1}{2} \right] 
\sim n^{2 - 2/\nu} \alpha_n^2 \frac{(2\nu)^{2/\nu}}{16\pi^2}.
\end{align*}
By assumption, $\alpha_n \sim \alpha n^{\frac{1}{\nu} - \frac{1}{2}}$ for some constant $\alpha \in (0, \infty)$. Thus, 
\begin{align*}
\lim_{n \rightarrow \infty} \frac{\mathds{E}[|X_n|^2]}{n} = \alpha^2 \frac{(2\nu)^{2/\nu}}{16\pi^2}.
\end{align*}

For the second moment of $|X_n|^2$, Parseval's theorem is applied again and gives that
\begin{align}\label{second}
\mathds{E}[(|X_n|^2)^2] &= 
\frac{1}{||K_n||_2^2} \int_{\R^n} (|x|^2 K_n(x))^2 {\mathrm d}x = \frac{1}{||K_n||_2^2} \int_{\R^n} \frac{(\triangle \hat{K}_n(\xi))^2}{(2\pi)^4} {\mathrm d}\xi.
\end{align}
Then, by the above computation of the Laplacian of $\hat{K}$, (\ref{important_integral}), and (\ref{powerexp_global}), 
\begin{align*}
&\mathds{E}[(|X_n|^2)^2] 
= \frac{\Gamma(\frac{n}{2} + 1) \alpha_n^n 2^{n/\nu}\nu^2 \alpha_n^{2\nu}}{(2\pi)^4\pi^{n/2}\Gamma(\frac{n}{\nu} + 1)}  \int_{\R^n} e^{(- 2| \alpha_n x |^{\nu})} \left(\nu \alpha_n^{\nu}| x|^{2\nu - 2} - ( \nu -2 + n) | x|^{\nu - 2}  \right)^2 {\mathrm d}x \\
& \qquad = \frac{\Gamma(\frac{n}{2} + 1) \alpha_n^n 2^{n/\nu}\nu^2 \alpha_n^{2\nu}}{(2\pi)^4\pi^{n/2}\Gamma(\frac{n}{\nu} + 1)} \bigg[(\nu \alpha_n^{\nu})^2 \int_{\R^n} e^{- 2| \alpha_n x |^{\nu}}| x|^{4\nu - 4} {\mathrm d}x \\
&\qquad \qquad - 2\nu \alpha_n^{\nu}( \nu -2 + n) \int_{\R^n}  e^{- 2| \alpha_n x |^{\nu}}| x|^{3\nu - 4} dx  + ( \nu -2 + n)^2 \int_{\R^n}   e^{- 2| \alpha_n x |^{\nu}} |x|^{2\nu - 4} {\mathrm d}x \bigg] \\
&\qquad = n\frac{\alpha_n^n2^{n/\nu}\nu^2 \alpha_n^{2\nu}}{(2\pi)^4\Gamma(\frac{n}{\nu} + 1)}  \bigg[\frac{(\nu \alpha_n^{\nu})^2 \Gamma(\frac{n + 4\nu - 4}{\nu})}{\nu 2^{(n + 4\nu - 4)/\nu}\alpha_n^{n + 4\nu - 4}} \\
&\qquad \qquad -  \frac{2\nu \alpha_n^{\nu}( \nu -2 + n)\Gamma(\frac{n + 3\nu - 4}{\nu})}{\nu 2^{(n + 3\nu - 4)/\nu}\alpha_n^{n + 3\nu - 4}}  + \frac{( \nu -2 + n)^2 \Gamma(\frac{n + 2\nu - 4}{\nu})}{\nu 2^{(n + 2\nu - 4)/\nu}\alpha_n^{n + 2\nu - 4}}\bigg] \\
&\qquad = \frac{n2^{4/\nu} \nu^2 \alpha_n^{4}}{(2\pi)^4\Gamma(\frac{n}{\nu} + 1)}  \bigg[ \frac{\nu\Gamma\left(\frac{n - 4}{\nu} + 4\right)}{2^{4}} - \frac{2( \nu -2 + n)\Gamma\left(\frac{n - 4}{\nu} + 3\right) }{2^3} +  \frac{( \nu -2 + n)^2\Gamma\left(\frac{n  - 4}{\nu} + 2\right) }{\nu 2^{2}} \bigg] \\
&\qquad = n\frac{2^{4/\nu} \alpha_n^4\Gamma\left(\frac{n  - 4}{\nu} + 1\right) }{(2\pi)^4\Gamma\left(\frac{n}{\nu} + 1\right)} \bigg[\frac{\nu^3}{2^{4}}\left(\frac{n  - 4}{\nu} + 3\right)\left(\frac{n  - 4}{\nu} + 2\right)\left(\frac{n  - 4}{\nu} + 1\right) \\ 
&\qquad \qquad  - \frac{\nu^2(n + \nu -2 )}{2^2}\left(\frac{n  - 4}{\nu} + 2\right)\left(\frac{n  - 4}{\nu} + 1\right) +   \frac{\nu( n + \nu -2 )^2}{2^{2}}\left(\frac{n  - 4}{\nu} + 1\right)  \bigg]  \\
&\qquad =   n\frac{2^{4/\nu} \alpha_n^4}{(2\pi)^4}\frac{ \Gamma\left(\frac{n  - 4}{\nu} + 1\right)}{\Gamma\left(\frac{n}{\nu} + 1\right)} \left(\frac{n^3}{2^4} - \frac{n^3}{2^2} + \frac{n^3}{2^2} + o(n^3) \right) \\
&\qquad =   n^4 \frac{2^{4/\nu}\alpha_n^4}{(2\pi)^4} \frac{ \Gamma\left(\frac{n  - 4}{\nu} + 1\right)}{\Gamma\left(\frac{n}{\nu} + 1\right)}  \left(\frac{1}{16} + o(1) \right) \\
&\qquad \sim n^4 \frac{2^{4/\nu}\alpha_n^4}{16(2\pi)^4}  \sqrt{\frac{\nu}{2\pi n}}\left(\frac{\nu e}{n}\right)^{\frac{n}{\nu}}  \sqrt{\frac{2\pi(n-4)}{\nu}} \left(\frac{n-4}{\nu e}\right)^{\frac{n-4}{\nu}} \\
&\qquad = n^4 \sqrt{\frac{n-4}{n}}\left(1-\frac{4}{n}\right)^{\frac{n}{\nu}} \left(\frac{n-4}{\nu e}\right)^{-\frac{4}{\nu}}\tfrac{\alpha_n^4 2^{4/\nu} }{16(2\pi)^4} 
\sim  n^4  \left(n-4\right)^{-\frac{4}{\nu}}\frac{\alpha_n^4 (2\nu)^{4/\nu}}{16(2\pi)^4}.
\end{align*}
Again, since $\alpha_n \sim \alpha n^{\frac{1}{\nu} - \frac{1}{2}}$, $\mathds{E}[(|X_n|^2)^2] = O(n^2)$, and
\begin{equation*}
\lim_{n\rightarrow \infty} \frac{\mathds{E}[(|X_n|^2)^2] }{n^{2}} =  \alpha^4 \frac{(2\nu)^{4/\nu}}{16(2\pi)^4}.
\end{equation*} 
Note that this limit is exactly the square of the limit of the expectation of $\frac{|X_n|^2}{n}$, implying
\begin{equation*}
\Var\left(\frac{|X_n|^2}{n^2}\right) =  \frac{\mathds{E}[(|X_n|^2)^2] }{n^{2}} - \left(\frac{\mathds{E}[|X_n|^2]}{n}\right)^2  \rightarrow 0 \text{ as } n \rightarrow \infty.
\end{equation*}
Thus, by Chebychev's inequality, $\frac{|X_n|}{\sqrt{n}} \rightarrow \alpha \frac{(2\nu)^{1/\nu}}{4\pi}$ in probability.

\section{Proof of Proposition \ref{bessel}}\label{bessel_app}

First, for $k \geq 0$, we see that
\begin{align*}
&\int_{\R^n} |x|^k K(x)^2 {\mathrm d}x = \int_{\R^n} |x|^k \left(e^{n\rho} 2^{(\sigma + n)/2} \Gamma\left( \frac{\sigma + n + 2}{2}\right) \frac{J_{(\sigma + n)/2}(2|x/\alpha|\sqrt{(\sigma + n)/2})}{(2|x/\alpha|\sqrt{(\sigma + n)/2})^{(\sigma + n)/2}}\right)^2 {\mathrm d}x \\
&\qquad = e^{2n\rho} 2^{(\sigma + n)} \Gamma\left( \frac{\sigma + n + 2}{2}\right)^2 \int_{\R^n} |x|^k  \frac{J_{(\sigma + n)/2}(2|x/\alpha|\sqrt{(\sigma + n)/2})^2}{(2|x/\alpha|\sqrt{(\sigma + n)/2})^{(\sigma + n)}} {\mathrm d}x \\
&\qquad= e^{2n\rho} 2^{(\sigma + n)} \Gamma\left( \frac{\sigma + n + 2}{2}\right)^2 \frac{2\pi^{n/2}}{\Gamma(\frac{n}{2})} \int_{0}^{\infty} r^{n-1} r^k  \frac{J_{(\sigma + n)/2}(2(r/\alpha)\sqrt{(\sigma + n)/2})^2}{(2(r/\alpha)\sqrt{(\sigma + n)/2})^{(\sigma + n)}} {\mathrm d}x,
\end{align*}
and by the change of variables $y = \left(\frac{2}{\alpha}\sqrt{\frac{\sigma + n}{2}}\right) r$, 
\begin{align*}
&= e^{2n\rho}  2^{\sigma + n}  \frac{2\pi^{n/2}\Gamma\left( \frac{\sigma + n + 2}{2}\right)^2}{\Gamma(\frac{n}{2})} \int_{0}^{\infty} \left(\frac{2}{\alpha}\sqrt{\frac{\sigma + n}{2}}\right)^{-k-n + 1} \frac{J_{(\sigma + n)/2}(y)^2}{y^{\sigma + 1-k}} \left(\frac{2}{\alpha}\sqrt{\frac{\sigma + n}{2}}\right)^{-1}{\mathrm d}y \\
&= e^{2n\rho} 2^{\sigma + n} \frac{2\pi^{n/2} \Gamma\left( \frac{\sigma + n + 2}{2}\right)^2\alpha^{k+n}}{\Gamma(\frac{n}{2})(2(\sigma + n))^\frac{k+n}{2}} \int_{0}^{\infty}\frac{J_{(\sigma + n)/2}(y)^2}{y^{\sigma +1 - k}} {\mathrm d}y.
\end{align*}
For $\sigma + 1 - k > 0$, from \cite[10.22.57]{DLMF},	
\begin{align*}
 \int_{0}^{\infty}\frac{J_{(\sigma + n)/2}(y)^2}{y^{\sigma +1 - k}} {\mathrm d}y 
= \frac{\Gamma\left(\frac{n}{2} + \frac{k}{2} \right) \Gamma(\sigma +1 - k)}{2^{\sigma - k + 1} \Gamma\left(\frac{\sigma - k}{2} + 1\right)^2 \Gamma\left(\sigma  - \frac{k}{2} + \frac{n}{2} + 1 \right)},
 \end{align*}
 and thus,
 \begin{align*}
\int_{\R^n} |x|^k K(x)^2 {\mathrm d}x &=e^{2n\rho} 2^{\sigma + n} \frac{2\pi^{n/2} \Gamma\left( \frac{\sigma + n + 2}{2}\right)^2\alpha^{k+n}}{\Gamma(\frac{n}{2})(2(\sigma + n))^\frac{k+n}{2}} \frac{\Gamma\left(\frac{n}{2} + \frac{k}{2} \right) \Gamma(\sigma +1 - k)}{2^{\sigma - k + 1} \Gamma\left(\frac{\sigma - k}{2} + 1\right)^2 \Gamma\left(\sigma  - \frac{k}{2} + \frac{n}{2} + 1 \right)} \\
  &= e^{2n\rho}\frac{(2\pi)^{n/2}\alpha^{k+n}2^{k/2} \Gamma\left( \frac{\sigma + n + 2}{2}\right)^2}{(\sigma + n)^{\frac{k +n}{2}} \Gamma(\frac{n}{2})} \frac{\Gamma\left(\frac{n}{2} + \frac{k}{2}\right) \Gamma(\sigma +1 - k)}{ \Gamma\left(\frac{\sigma - k}{2} + 1\right)^2 \Gamma\left(\sigma  - \frac{k}{2} + \frac{n}{2} + 1 \right)}.
 \end{align*}
 Then, for $\sigma > 0$,
 \begin{align*}
 \mathds{E}[|X_n|] &= \tfrac{1}{||K_n||_2^2} \int_{\R^n} |x| K_n(x)^2 {\mathrm d}x \\
 &= \frac{(2\pi)^{\frac{n}{2}}\alpha^{1+n}2^{1/2}\Gamma\left( \frac{\sigma + n + 2}{2}\right)^2\Gamma\left(\frac{n}{2} + \frac{1}{2}\right) \Gamma(\sigma)}{(\sigma + n)^{\frac{1+n}{2}} \Gamma(\frac{n}{2}) \Gamma\left(\frac{\sigma + 1}{2} \right)^2 \Gamma\left(\sigma  - \frac{1}{2} + \frac{n}{2} + 1 \right)} \frac{(\sigma + n)^{\frac{n}{2}}  \Gamma\left(\frac{\sigma}{2} + 1\right)^2 \Gamma\left(\sigma + \frac{n}{2} + 1 \right)}{(2\pi)^{\frac{n}{2}}\alpha^{n}\Gamma(\sigma +1)\Gamma\left( \frac{\sigma}{2} + \frac{n}{2} + 1 \right)^2} \\
&=  \frac{\alpha 2^{1/2}}{(\sigma + n)^{1/2} \Gamma(\frac{n}{2})} \frac{\Gamma\left(\frac{n}{2} + \frac{1}{2}\right) \Gamma(\sigma )\Gamma\left(\frac{\sigma}{2} + 1\right)^2 \Gamma\left(\sigma + \frac{n}{2} + 1 \right)}{ \Gamma\left(\frac{\sigma}{2} + \frac{1}{2}\right)^2 \Gamma\left(\sigma + \frac{n}{2} + \frac{1}{2} \right)\Gamma(\sigma +1)} \\
&\sim  \frac{\alpha 2^{1/2}}{(\sigma + n)^{1/2} \Gamma(\frac{n}{2})} \frac{\Gamma\left(\frac{n}{2}\right)\left(\frac{n}{2}\right)^{\frac{1}{2}} \Gamma(\sigma)\Gamma\left(\frac{\sigma}{2} + 1\right)^2 \Gamma\left(\frac{n}{2}\right)\left(\frac{n}{2}\right)^{\sigma + 1}}{ \Gamma\left(\frac{\sigma +1}{2}\right)^2 \Gamma\left(\frac{n}{2}\right)\left(\frac{n}{2}\right)^{\sigma  + \frac{1}{2} } \Gamma(\sigma +1)} \\
&=  \frac{\alpha 2^{1/2}}{(\sigma + n)^{1/2} } \frac{\left(\frac{n}{2}\right) \Gamma(\sigma )\Gamma\left(\frac{\sigma}{2} + 1\right)^2}{ \Gamma\left(\frac{\sigma +1}{2} \right)^2 \Gamma(\sigma +1)} 
\sim  n^{1/2} \frac{\alpha}{2^{1/2}} \frac{\Gamma(\sigma)\Gamma\left(\frac{\sigma}{2} + 1\right)^2}{ \Gamma\left(\frac{\sigma +1}{2}\right)^2 \Gamma(\sigma +1)} = O(n^{\frac{1}{2}}).
\end{align*}
Now, let $\beta > \frac{1}{2}$. By Markov's inequality, 
\begin{align*}
\lim_{n \to \infty} \frac{\mathds{E}[\eta_n(B_n(Rn^{\beta})^c)]}{\mathds{E}[\eta_n(\R^n)]} &= \lim_{n \to \infty} \mathds{P}\left(|X_n| \geq R n^{\beta} \right) \leq \lim_{n \to \infty} \frac{\mathds{E}|X_n|}{Rn^{\beta}} = 0.
\end{align*} 

\added{\section{Proof of Proposition \ref{Whittle}}

First, from \cite[6.576.3]{ryzhik}, we have for all $\nu > 0$ and $k > 2\nu - 1$, 
\begin{align}\label{e:bess2_int}
\int_0^{\infty} r^k \mathbb{K}_{\nu}\left(\frac{r}{\alpha}\right)^2 dr = \frac{2^{-2+k}\alpha^{k + 1}}{\Gamma(1 + k)}\Gamma\left(\frac{1 + k}{2} + \nu\right)\Gamma \left(\frac{k + 1}{2}\right)^2\Gamma\left(\frac{1 + k}{2} - \nu\right),
\end{align}
where $\mathbb{K}_{\nu}$ is the modified Bessel function of the second kind.

For the Whittle-Mat\'{e}rn Kernel \eqref{e:whittle_K},
\begin{align*}
\int_{\R^n} K_n(x)^2 dx &= \int_{\R^n} e^{2n\rho} \frac{2^{2 - 2\nu}}{\Gamma(\nu)^2} \frac{|x|^{2\nu}}{\alpha^{2\nu}}\mathbb{K}_{\nu}\left(\frac{|x|}{\alpha}\right)^2 dx \\
&= \frac{2\pi^{\frac{n}{2}}}{\Gamma(\frac{n}{2})}e^{2n\rho} \frac{2^{2 - 2\nu}}{\Gamma(\nu)^2 \alpha^{2\nu}} \int_0^{\infty} r^{n-1}r^{2\nu} \mathbb{K}_{\nu}(r)^2 dr \\
&= e^{2n\rho}\frac{2\pi^{\frac{n}{2}}}{\Gamma(\frac{n}{2})}\frac{2^{2 - 2\nu}}{\Gamma(\nu)^2 \alpha^{2\nu}} \int_0^{\infty} r^{n-1 + 2\nu} \mathbb{K}_{\nu}(r)^2 dr.
\end{align*}
Then by \eqref{e:bess2_int},
\begin{align*}
 \int_0^{\infty} r^{n-1 + 2\nu} \mathbb{K}_{\nu}(r)^2 dr &=  \frac{2^{-3+n + 2 \nu }\alpha^{n + 2\nu}}{\Gamma(n + 2\nu)}\Gamma\left(\frac{n + 2\nu}{2} + \nu\right)\Gamma \left(\frac{n + 2\nu}{2}\right)^2\Gamma\left(\frac{n + 2\nu}{2} - \nu\right) \\
 &=  \frac{2^{-3+n + 2 \nu }\alpha^{n + 2\nu}}{\Gamma(n + 2\nu)}\Gamma\left(\frac{n}{2} + 2\nu\right)\Gamma \left(\frac{n}{2} + \nu \right)^2\Gamma\left(\frac{n}{2}\right).
\end{align*}
Similarly,
\begin{align*}
\int_{\R^n} |x|^2 K_n(x)^2 dx &= e^{2n\rho}\frac{2\pi^{\frac{n}{2}}}{\Gamma(\frac{n}{2})}\frac{2^{2 - 2\nu}}{\Gamma(\nu)^2 \alpha^{2\nu}} \int_0^{\infty} r^{n+1 + 2\nu} \mathbb{K}_{\nu}(r)^2 dr.
\end{align*}
and also by \eqref{e:bess2_int},
\begin{align*}
 \int_0^{\infty} r^{n+1 + 2\nu} \mathbb{K}_n(r)^2 dr = \frac{2^{-1+n + 2 \nu }\alpha^{n + 2 + 2\nu}}{\Gamma(n + 2 + 2\nu)}\Gamma\left(\frac{n}{2} + 2\nu + 1\right)\Gamma \left(\frac{n}{2} + \nu + 1\right)^2\Gamma\left(\frac{n}{2} + 1\right)
\end{align*}
Then,
\begin{align*}
\mathds{E}[|X_n|^2] &= \frac{\int_{\R^n} |x|^2 K_n(x)^2 dx}{\int_{\R^n} K_n(x)^2 dx} 
= \frac{ (2\alpha)^{2}\Gamma(n + 2\nu)\Gamma\left(\frac{n}{2} + 2\nu + 1\right)\Gamma \left(\frac{n}{2} + \nu + 1\right)^2\Gamma\left(\frac{n}{2} + 1\right)}{\Gamma(n + 2 + 2\nu)\Gamma\left(\frac{n}{2} + 2\nu\right)\Gamma \left(\frac{n}{2} + \nu \right)^2\Gamma\left(\frac{n}{2}\right)} \\
&= \frac{ (2\alpha)^{2}\left(\frac{n}{2} + 2\nu\right) \left(\frac{n}{2} + \nu\right)^2\left(\frac{n}{2} \right)}{(n + 1 + 2\nu)(n + 2\nu)} \sim \left(\frac{\alpha}{2}\right)^2n,
\end{align*}
as $n \rightarrow \infty$, and this implies
\[\frac{\mathds{E}[|X_n|^2]^{\frac{1}{2}}}{\sqrt{n}} \rightarrow \frac{\alpha}{2}, \text{ as } n \rightarrow \infty.\]

Thus, since the Whittle Mat\'{e}rn kernel is log-concave, the conclusion holds by Theorem \ref{logconcave}.}

\section{Proof of Proposition \ref{Cauchy}}

First, recall the the beta function satisfies
\begin{align*}
B(x,y) := \int_0^1 t^{x-1}(1 - t)^{y-1} {\mathrm d}t = \int_0^{\infty} t^{x-1} (1 + t)^{-(x+y)} {\mathrm d}t = \frac{\Gamma(x)\Gamma(y)}{\Gamma(x + y)}.
\end{align*}
Then, for any $k \geq 0$,
\begin{align*}
\int_{\R^n} |x|^k K_n(x)^2 {\mathrm d}x &= \int_{\R^n} |x|^k\frac{e^{2n\rho}}{(1 + |\frac{x}{\alpha_n}|^2)^{2\nu + n}} {\mathrm d}x \\
&= e^{2n\rho} \frac{2\pi^{n/2}}{\Gamma(\frac{n}{2})} \int_0^{\infty} r^{n-1 + k}\left(1 + \frac{r^2}{\alpha_n^2}\right)^{-2\nu - n} {\mathrm d}r \\
&= e^{2n\rho} \frac{\pi^{n/2}}{\Gamma(\frac{n}{2})} \alpha_n^{n+k} \int_0^{\infty} t^{\frac{n}{2}-1 + \frac{k}{2}}(1 + t)^{-(2\nu + n)}  {\mathrm d}t \\
&= e^{2n\rho} \frac{\pi^{n/2}}{\Gamma(\frac{n}{2})} \alpha_n^{n+k} B\left(\frac{n}{2} + \frac{k}{2}, 2\nu + \frac{n}{2} - \frac{k}{2}\right).
\end{align*}
Thus, the expectation of $|X_n|^2$ is
\begin{align*}
\mathds{E}[|X_n|^2] &= \frac{1}{||K_n||_2^2} \int_{\R^n} |x|^2 K_n(x)^2 {\mathrm d}x
= \alpha_n^2 \frac{B(\frac{n}{2} + 1, 2\nu + \frac{n}{2} - 1)}{B(\frac{n}{2}, 2\nu + \frac{n}{2})} \\
&= \alpha_n^2 \frac{\Gamma(\frac{n}{2} + 1)\Gamma(2\nu + \frac{n}{2} - 1)\Gamma(n + 2\nu)}{\Gamma(n + 2\nu) \Gamma(\frac{n}{2})\Gamma(2\nu + \frac{n}{2})} = \alpha_n^2 \frac{n}{2(\frac{n}{2} + 2\nu - 1)} = \alpha_n^2 \frac{n}{n + 4\nu - 2},
\end{align*}
and 
\begin{align*}
\mathds{E}[|X_n|^4] 
&= \alpha_n^4 \frac{B(\frac{n}{2} + 2, 2\nu +\frac{n}{2} - 2)}{B(\frac{n}{2}, 2\nu + \frac{n}{2})} 
= \alpha_n^4\frac{\Gamma(\frac{n}{2} + 2)\Gamma(2\nu + \frac{n}{2} - 2)\Gamma(n + 2\nu)}{\Gamma(n + 2\nu) \Gamma(\frac{n}{2})\Gamma(2\nu + \frac{n}{2})} \\
&= \alpha_n^4 \frac{(\frac{n}{2} + 1)\frac{n}{2}}{(2\nu + \frac{n}{2} -2)(2\nu + \frac{n}{2} - 1)} 
= \alpha_n^4 \frac{n(n + 2)}{(n + 4\nu -4)(n + 4\nu - 2)}.
\end{align*}
Thus, by the assumption that $\alpha_n \sim \alpha n^{\frac{1}{2}}$ as $n \rightarrow \infty$ for some $\alpha > 0$,
\begin{align*}
\lim_{n \rightarrow \infty}\frac{\mathds{E}[|X_n|^2]}{n} = \alpha^2 \text{ and }  \lim_{n \rightarrow \infty} \frac{\Var(|X_n|^2)}{n^2}  = 0.
\end{align*}
Thus, by Chebychev's inequality, $\frac{|X_n|}{\sqrt{n}} \rightarrow \alpha$ in probability.

\section{Proof of Proposition \ref{gauss_beta2}}

By Proposition \ref{Laguerre_rate},
\begin{align*}
\lim_{n \rightarrow \infty} - \frac{1}{n} \ln \mathds{E}[\eta_n(B_n(\sqrt{n}R))] = \begin{cases} - \rho - \frac{1}{2} \log 2\pi e  + \frac{2R^2}{\alpha^2m} - \log R, & 0 < R < \sqrt{m}\frac{\alpha}{2} \\ - \rho - \log \alpha - \frac{1}{2} \log \frac{m\pi}{2} , & R > \sqrt{m}\frac{\alpha}{2}.
\end{cases}
\end{align*}
Recall that $\lim_{n \rightarrow \infty}  \frac{1}{n} \ln \mathds{E}[\Phi_n(B_n(\sqrt{n}R))] = \rho + \frac{1}{2} \log 2 \pi e + \log R$. Thus,
\begin{align*}
&\lim_{n \rightarrow \infty} - \frac{1}{n} \ln  \frac{\mathds{E}[\eta_n(B_n(\sqrt{n}R))]}{\mathds{E}[\Phi_n(B_n(\sqrt{n}R))] }  \\
&= \lim_{n \rightarrow \infty} - \frac{1}{n} \ln \mathds{E}[\eta_n(B_n(\sqrt{n}R))] + \frac{1}{n} \ln \mathds{E}[\Phi_n(B_n(\sqrt{n}R))] \\
&= \begin{cases} - \rho - \frac{1}{2} \log 2\pi e  + \frac{2R^2}{\alpha^2m} - \log R + \rho + \frac{1}{2}\log 2\pi e + \log R, & 0 < R < \sqrt{m}\frac{\alpha}{2} \\ - \rho - \log \alpha - \frac{1}{2} \log \frac{m\pi}{2} +\rho + \frac{1}{2}\log 2\pi e + \log R, & R > \sqrt{m}\frac{\alpha}{2} \end{cases} \\
&= \begin{cases}  \frac{2R^2}{\alpha^2m}, & 0 < R < \sqrt{m}\frac{\alpha}{2} \\  \frac{1}{2} +\log 2 - \log \alpha - \frac{1}{2} \log m + \log R, & R > \sqrt{m}\frac{\alpha}{2}.
\end{cases}
\end{align*}

\acks
The work of both authors was supported by a grant of the Simons Foundation (\#197982 to UT Austin). The work of the second author was supported by the National Science Foundation Graduate Research Fellowship under Grant No. DGE-1110007.

\bibliography{RefDPP}

\begin{thebibliography}{10}

\bibitem{Venkat}
{\sc Anantharam, V. and Baccelli, F.} (2016).
\newblock The {Boolean} model in the {Shannon} regime: Three thresholds and
  related asymptotics.
\newblock {\em Journal of Applied Probability\/} {\bf 53,} 1001--1018.

\bibitem{Lav}
{\sc {Biscio}, C. and {Lavancier}, F.} (2016).
\newblock {Quantifying repulsiveness of determinantal point processes}.
\newblock {\em Bernoulli\/} {\bf 22,} 2001--2028.

\bibitem{Blasz}
{\sc Blaszczyszyn, B. and Yogeshwaran, D.} (2014).
\newblock On comparison of clustering properties of point processes.
\newblock {\em Advances in Applied Probability\/} {\bf 46,} 1--20.

\bibitem{Stoyan}
{\sc Chiu, S.~N., Stoyan, D., Kendall, W.~S. and Mecke, J.} (2013).
\newblock {\em Stochastic Geometry and its Applications} third~ed.
\newblock Wiley.

\bibitem{LDP}
{\sc Dembo, A. and Zeitouni, O.} (1998).
\newblock {\em Large Deviations Techniques and Applications} 2nd~ed.
\newblock Springer.

\bibitem{Guedon}
{\sc Fradelizi, M., Gu\'{e}don, O. and Pajor, A.} (2014).
\newblock Thin-shell concentration for convex measures.
\newblock {\em Studia Mathematica\/} {\bf 223,}.

\bibitem{Goldman}
{\sc Goldman, A.} (2010).
\newblock The {Palm} measure and the {Voronoi} tessellation for the {Ginibre}
  process.
\newblock {\em Annals of Applied Probability\/} {\bf 20,} 90--128.

\bibitem{ryzhik}
{\sc Gradshteyn, I. and Ryzhik, I.} (2007).
\newblock {\em Table of Integrals, series, and products} seventh~ed.
\newblock Elsevier.

\bibitem{Milman}
{\sc Gu{\'e}don, O. and Milman, E.} (2011).
\newblock Interpolating thin-shell and sharp large-deviation estimates for
  isotropic log-concave measures.
\newblock {\em Geometric and Functional Analysis\/}.

\bibitem{Peres}
{\sc Hough, J.~B., Krishnapur, M., Peres, Y. and Virag, B.} (2009).
\newblock Zeros of {Gaussian} analytic functions and determinantal point
  processes.
\newblock {\em American Mathematical Society\/} {\bf 51,}.

\bibitem{Klartag}
{\sc Klartag, B.} (2007).
\newblock A central limit theorem for convex sets.
\newblock {\em Inventiones Mathematicae\/} {\bf 168,} 91--131.

\bibitem{Taskar}
{\sc Kulesza, A. and Taskar, B.} (2012).
\newblock Determinantal point processes for machine learning.
\newblock {\em Foundations and Trends in Machine Learning\/} {\bf 5,} 123--286.

\bibitem{Kuna}
{\sc Kuna, T., Lebowitz, J. and Speer, E.} (2007).
\newblock Realizability of point processes.
\newblock {\em Journal of Statistical Physics\/} {\bf 129,} 417--439.

\bibitem{Moller_Arxiv}
{\sc Lavancier, F., M{\o}ller, J. and Rubak, E.} (2012).
\newblock Determinantal point process models and statistical inference:
  Extended version.
\newblock {\em arXiv:1205.4818\/}.

\bibitem{Moller}
{\sc Lavancier, F., M{\o}ller, J. and Rubak, E.} (2015).
\newblock Determinantal point process models and statistical inference.
\newblock {\em Journal of the Royal Statistical Society: Series B (Statistical
  Methodology)\/} {\bf 77,} 853--877.

\bibitem{Baccelli_wireless}
{\sc Li, Y., Baccelli, F., Dhillon, H.~S. and Andrews, J.~G.} (2015).
\newblock Statistical modeling and probabilistic analysis of cellular networks
  with determinantal point processes.
\newblock {\em IEEE Transactions on Communications\/} {\bf 63,} 3405--3422.

\bibitem{Macchi}
{\sc Macchi, O.} (1975).
\newblock The coincidence approach to stochastic point processes.
\newblock {\em Advances in Applied Probability\/} {\bf 7,} 83--122.

\bibitem{JMEO}
{\sc M{\o}ller, J. and O'Reilly, E.} (2018).
\newblock Couplings for determinantal point processes and their reduced {Palm}
  distribution with a view to quantifying repulsiveness.
\newblock {\em arXiv:1806.07347\/}.

\bibitem{DLMF}
{NIST Digital Library of Mathematical Functions}.
\newblock http://dlmf.nist.gov/, Release 1.0.14 of 2016-12-21.
\newblock F.~W.~J. Olver, A.~B. {Olde Daalhuis}, D.~W. Lozier, B.~I. Schneider,
  R.~F. Boisvert, C.~W. Clark, B.~R. Miller and B.~V. Saunders, eds.

\bibitem{Ripley}
{\sc Ripley, B.~D.} (1976).
\newblock The second-order analysis of stationary point process.
\newblock {\em Journal of Applied Probability\/} {\bf 13,} 255--266.

\bibitem{Rudin}
{\sc Rudin, W.} (1991).
\newblock {\em Functional Analysis} second~ed.
\newblock International Series in Pure and Applied Mathematics. McGraw-Hill,
  Inc.

\bibitem{Weil}
{\sc Schneider, R. and Wolfgang, W.} (2008).
\newblock {\em Stochastic and Integral Geometry}.
\newblock Probability and Its Applications. Springer.

\bibitem{Shirai}
{\sc Shirai, T. and Takahashi, Y.} (2003).
\newblock Random point fields associated with certain {Fredholm} determinants
  {I}: {Fermion}, {Poisson}, and {Boson} point processes.
\newblock {\em Journal of Functional Analysis\/} {\bf 205,} 414--463.

\bibitem{Torquato}
{\sc Torquato, S., Scardicchio, A. and Zachary, C.~E.} (2008).
\newblock Point processes in arbitrary dimension from {Fermionic} gases, random
  matrix theory, and number theory.
\newblock {\em Journal of Statistical Mechanics\/}.

\bibitem{Yu}
{\sc Yu, Y.} (2017).
\newblock On normal variance-mean mixtures.
\newblock {\em Statistics and Probability Letters\/} {\bf 121,} 45--50.

\end{thebibliography}
\bibliographystyle{apt}

\end{document}